\documentclass{amsart}
\usepackage{amssymb}
\usepackage[latin1]{inputenc}
\usepackage[swedish,english]{babel}
\usepackage{ifpdf}
\ifpdf
  \usepackage[pdftex]{graphicx}
\else
  \usepackage[dvips]{graphicx}
\fi
\usepackage{subfigure}
\usepackage[pdftitle={Surfactants in two-phase fluid flow},
  pdfauthor={Stefan Engblom},
  pdffitwindow=true,
  breaklinks=true,
  colorlinks=true,
  urlcolor=blue,
  linkcolor=red,
  citecolor=red,
  anchorcolor=red]{hyperref}
\usepackage[numbers]{natbib}

\numberwithin{equation}{section}
\numberwithin{table}{section}
\numberwithin{figure}{section}

\theoremstyle{plain}

\theoremstyle{definition}

\theoremstyle{remark}
\newtheorem*{conclusion}{Conclusion}

\newtheorem*{example}{Example}

\newcommand{\fatu}{\mathbf{u}}

\newcommand{\Fphi}{F_{\phi}}
\newcommand{\Fpsi}{F_{\psi}}
\newcommand{\Fex}{F_{\mbox{\tiny{ex}}}}
\newcommand{\Fone}{F_{1}}

\newcommand{\Cn}{\operatorname{Cn}}
\newcommand{\Ex}{\operatorname{Ex}}
\newcommand{\Pii}{\operatorname{Pi}}
\newcommand{\Pe}{\operatorname{Pe}}
\newcommand{\Reynold}{\operatorname{Re}}
\newcommand{\Ca}{\operatorname{Ca}}

\newcommand{\Ordo}[1]{\mathcal{O}\left( #1 \right)}
\def\sech{\qopname\relax o{sech}} 

\newcommand{\Hone}{H^{1}}
\newcommand{\Ltwo}{L^{2}}
\newcommand{\Leg}{P}
\newcommand{\phihat}{\hat{\phi}}
\newcommand{\Phihat}{\hat{\Phi}}
\newcommand{\psihat}{\hat{\psi}}

\newcommand{\TOL}{\operatorname{TOL}}
\newcommand{\Atol}{\operatorname{Atol}}
\newcommand{\Rtol}{\operatorname{Rtol}}
\newcommand{\err}{\operatorname{err}}
\newcommand{\cerr}{\operatorname{cerr}}

\newcommand{\eq}{\mbox{{\small eq}}}
\newcommand{\eeq}{\mbox{{\tiny eq}}}
\newcommand{\psiin}{\psi_{\mbox{{\tiny in}}}}
\newcommand{\psiout}{\psi_{\mbox{{\tiny out}}}}

\begin{document}

\title[Surfactants in two-phase fluid flow]{On diffuse interface
  modeling and simulation of surfactants in two-phase fluid flow}

\author[S. Engblom]{Stefan Engblom}
\address[S. Engblom]{Division of Scientific Computing \\
  Department of Information Technology \\
  Uppsala University \\
  SE-751 05 Uppsala, Sweden.}

\urladdr[S. Engblom]{\url{http://user.it.uu.se/~stefane}}
\email{stefane@it.uu.se}

\author[M. Do-Quang]{Minh Do-Quang}
\address[M. Do-Quang \and G. Amberg]{Linn{\'e} Flow Centre \\
  Department of Mechanics \\
  School of Engineering Science \\
  Royal Institute of Technology \\
  S-100 44 Stockholm \\
  Sweden}
\email{minh@mech.kth.se, gustava@mech.kth.se}

\author[G. Amberg]{Gustav Amberg}

\author[A-K. Tornberg]{Anna-Karin Tornberg}
\address[A-K. Tornberg]{Linn{\'e} Flow Centre \\
  Department of Numerical Analysis \\
  School of Computer Science and Communication \\
  Royal Institute of Technology \\
  S-100 44 Stockholm \\
  Sweden}
\email{annak@nada.kth.se}

\thanks{Corresponding author: S. Engblom, telephone +46-18-471 27 54,
  fax +46-18-51 19 25.}

\date{October 5, 2012}

%
%
%
%


\begin{abstract}

  An existing phase-field model of two immiscible fluids with a single
  soluble surfactant present is discussed in detail. We analyze the
  well-posedness of the model and provide strong evidence that it is
  mathematically ill-posed for a large set of physically relevant
  parameters. As a consequence, critical modifications to the model
  are suggested that substantially increase the domain of
  validity. Carefully designed numerical simulations offer informative
  demonstrations as to the sharpness of our theoretical results and
  the qualities of the physical model. A fully coupled hydrodynamic
  test-case demonstrates the potential to capture also non-trivial
  effects on the overall flow.

\end{abstract}

\selectlanguage{english}

\maketitle


\section{Introduction}
\label{sec:intro}

The presence of surface active substances may greatly affect the
physical properties of fluid mixtures. Indeed, these effects are used
critically in many important applications in everyday life; detergents
and oil-water emulsions in food are two immediate examples. The fact
that surfactants \emph{lower} the surface tension is exploited in both
of these cases: detergents make the water more ``wet'', and an
emulsifying agent stabilizes an emulsion by preventing small droplets
to coalesce.

For humans, probably the most critical everyday usage of surfactants
is made in the \textit{alveoli} in the lung, where \textit{pulmonary
  surfactant}, amongst other things, prevents lung collapse at the end
of expiration \cite{surflungs,surfalveol}.

An interesting and very striking realization \textit{in vivo} was
reported recently in \cite{maze} where \emph{chemotaxis} was
implemented for small droplets of fluid in a bulk solution, physically
contained in a maze. The net transport at the millimeter-scale, and
the subsequent solution to the maze, was achieved by a clever usage of
surfactant and a pre-existing pH-gradient. Quite likely, such a
constructive set-up could find applications in lab-on-a-chip
manufacturing.

Consisting of hydrophobic ``heads'' and hydrophilic ``tails'',
surfactant molecules have a strong preference to occupy sites at the
water-fluid or water-gas interfaces. Below the \emph{critical micelle
  concentration} (CMC), surfactants therefore adsorb efficiently to
the interfaces where their physical effects become prominent. Above
the CMC, additionally, spontaneous formation of stable groups of
surfactants --- \emph{micelles} --- occurs in the bulk solution
\cite{micelle}.

Given the `thermodynamical' signature of surfactants; that is, the
diffusion-limited flow into the interfaces and the spontaneous
creation of highly regular micelles from an unordered state in the
bulk, modeling through some kind of system's energy assumption is a
tempting approach.

Ariel, Diamant, and Andelman \cite{surfkinetics,surfkinetics0} have
successfully postulated free-energy terms with theoretically
convincing properties. Their approach is inherently \emph{sharp} in
that they set up equations for the \emph{interface}, the
\emph{sub-surface region}, and the \emph{bulk}. A later work
\cite{micelle} shows that this methodology can be extended to include
also the region above the CMC. A great feature with this type of
modeling is the fact that all properties of the model result from a
single postulated entity; the system's free energy.

Given the multitude of scales present and the complex coupling to
hydrodynamics, numerical experiments become important as tools to gain
a better understanding and a fuller physical insight.

In a so-called sharp interface method, the interface is considered to
be infinitesimally thin, and its exact location is represented either
explicitly (front-tracking with e.g.~Lagrangian markers), or
implicitly (e.g.~level-set \cite{levelset_surf} and volume of fluid
(VOF) methods \cite{surfVOF}). The evolution of the surfactant
concentration on each interface can be described by a partial
differential equation on this time-dependent manifold
\cite{sharp_surfactant}. Techniques for solving this PDE and include
insoluble surfactants in multiphase flow simulations have been
developed based on several different interface representation
techniques \cite{drop_deformation,surfVOF,
  insoluble_surf,immersed_surf, marangoni_drop,levelset_surf}.

When the surfactants are soluble also in the bulk, source terms due to
adsorption and desorption terms enter this PDE, and this must be
coupled to a PDE for the bulk concentration of surfactants with
appropriate boundary conditions for the surfactant flux. See
e.g.~Eggleton and Stebe \cite{eggleton_surf} for an early
reference. To simplify matters for simulations, it is in this
reference assumed that the adsorption and desorption from the bulk is
diffusion dominated and that the bulk surfactant concentration is
spatially constant. Work to consider the full problem has started only
recently and is currently an active area of research \cite{surfhybrid,
  Khatri_thesis,front_tracking_surf2,surfsoluble,front_tracking_surf}.

Given a thermodynamic description as a postulated free energy of the
system, including surface effects introduced via gradient energy terms
\cite{ternaryfluid,surfdynamics,diffsurf2,diffsurf,
  dropletsurf,diffemul}, it is attractive to represent the interface
as a rapid but continuous transition in a concentration variable. This
is the essence of \emph{phase-field} or \emph{diffuse interface
  modeling}. Without the need to explicitly track the interface, even
complicated topological transformations can be studied, and the
coupling to the full hydrodynamic flow is fairly straightforward
\cite{CahnHilliard}.

However, this has only recently begun to be explored for surfactant
laden interfaces. In \cite{surfsoluble}, a phase-field method was used
to represent the interfaces. The surfactant treatment was introduced
in a manner quite natural to a sharp interface method, and was not
included in a thermodynamical derivation of a diffuse interface method
for the full problem. In \cite{diffsurf}, and with some additions, in
\cite{diffsurf2}, a seemingly natural diffuse interface model for
surfactant laden liquid-liquid mixtures is presented, and solved with
lattice Boltzmann methods. The model is verified against classical
results, and some cases of convection dominated flow, such as a
deforming droplet in a shear flow are presented. This formulation is
promising from a computational and physical point of view, but it is
rather complex and its mathematical properties have not been
investigated. We will show below that this model can in fact be
expected to be mathematically ill-posed under certain conditions. The
purpose of the present paper is to clarify the conditions under which
the model is ill-posed, and also to present variants that behave
better.

In Section~\ref{sec:analysis} we thus present the analysis of the
model, and we also derive three alternatives in the form of
modifications to the original model. The properties of these
alternatives are investigated in Section~\ref{sec:numexp1} using a
highly accurate one-dimensional numerical scheme, and in a more
qualitative sense in two spatial dimensions in
Section~\ref{sec:numexp2}. An outlook with conclusions is finally
found in Section~\ref{sec:conclusions}.


\section{Diffuse interface model for surfactant flow}
\label{sec:model}

In this section we shall discuss our `baseline' phase-field model for
multiphase flow incorporating surfactants; --- later we shall have
reasons to alter the actual model by incorporating new terms. Our
starting point is a model originally presented in \cite{diffsurf}, but
incorporating the work also of others.

Specifically, a diffuse interface model for surfactant-controlled
emulsification was proposed in \cite{diffemul} and a related model for
oil-water-detergent mixtures appeared in \cite{ternaryfluid}. A model
with more easily understood interface adsorption properties was thus
proposed in \cite{diffsurf}, incorporating free energy terms developed
previously for a sharp model \cite{sharpsurf}, but drawing also on
some earlier work \cite{surfdynamics}. For clarity we take a
non-dimensional version as our baseline model (see
Appendix~\ref{app:nondim} for the precise relation to the original
model in \cite{diffsurf}).

The degrees of freedom are the \emph{phase-field variable} $\phi \in
[-1,1]$, the \emph{surfactant volume fraction} $\psi \in [0,1]$, and
the \emph{fluid velocity field} $\fatu$. The governing equations for
$\phi$ and $\psi$ take the form of Cahn-Hilliard-type equations
\cite{CahnHilliard},
\begin{align}
  \label{eq:phiPDE}
  \frac{\partial \phi}{\partial t}+\nabla \cdot (\phi \fatu) &=
  \frac{1}{\Pe_{\phi}} \nabla \cdot M_{\phi} \nabla \mu_{\phi} =
  \frac{1}{\Pe_{\phi}} \Delta \mu_{\phi}, \\
  \label{eq:psiPDE}
  \frac{\partial \psi}{\partial t}+\nabla \cdot (\psi \fatu) &=
  \frac{1}{\Pe_{\psi}}
  \nabla \cdot M_{\psi} \nabla \mu_{\psi},
\end{align}
in terms of \emph{P\'eclet} numbers $\Pe$, \emph{mobilities} $M$, and
\emph{chemical potentials} $\mu$ which will be prescribed below. As
indicated in \eqref{eq:phiPDE} the model is somewhat simplified in
that it assumes the ``shallow quench'' limit \cite[Sect.~9]{Ch4CHE},
and hence that the mobility $M_{\phi}$ for $\phi$ is assumed to be
constant (and in fact normalized to unity by scaling the P\'eclet
number $\Pe_{\phi}$ appropriately). This simplification is quite
standard and is adopted here mainly for convenience.

Eqs.~\eqref{eq:phiPDE}--\eqref{eq:psiPDE} are coupled to the
Navier-Stokes equations in the form
\begin{align}
  \label{eq:cont}
  \nabla \cdot \fatu &= 0, \\
  \label{eq:NS}
  \rho \left( \frac{\partial \fatu}{\partial t}+
    \fatu \cdot \nabla \fatu \right) &= -\nabla P+
  \frac{1}{\Reynold}\nabla \cdot 
  \left(\rho \nu [\nabla \fatu+(\nabla \fatu)^{T}] \right) \\
  \nonumber
  &\phantom{=} -\frac{1}{\Ca\Cn\Reynold}
  \left( \phi\nabla\mu_{\phi}+\psi\nabla\mu_{\psi} \right)
\end{align}
where the pressure tensor $P$ enforces the incompressibility condition
\eqref{eq:cont}, and where $\Reynold$ is the Reynolds number, $\Ca$
the Capillary number, and where, respectively, $\rho$ and $\nu$ denote
the fluid's density and kinematic viscosity. In the general case,
$\rho$ and $\nu$ may depend on $\phi$ as well as on other state
variables, requiring that additional effects be taken into account
\cite{densityCH}. For our present purposes and for simplicity they are
assumed to be constants.

\subsection{The Ginzburg-Landau free energy}

Recall that in diffuse models of binary fluids in equilibrium, the
phase-field variable $\phi$ typically describes a planar interface at
the origin of the form
\begin{align}
  \phi(x) &= \tanh (x/\Cn),
\end{align}
where $\Cn$ is the \emph{Cahn} number expressing the ratio between the
interface width and the characteristic length scale.

In \eqref{eq:phiPDE}--\eqref{eq:psiPDE} the chemical potentials are
derived from a Ginzburg-Landau \emph{free energy functional},
\begin{align}
  \label{eq:F}
  \int_{\Omega} F \, dx &\equiv \int_{\Omega}
  \Fphi+\Fpsi+\Fone+\Fex \, dx.
\end{align}
The classical constant mobility Cahn-Hilliard potential
\cite[Sect.~4]{Ch4CHE} is given by
\begin{align}
  \label{eq:Fphi}
  \Fphi &= -\frac{\phi^{2}}{2}+\frac{\phi^{4}}{4}+
  \frac{\Cn^{2}}{4}(\nabla \phi)^{2},
\end{align}
where we note that the polynomial part is a ``double well'' potential,
$[(\phi^{2}-1)^{2}-1]/4$, and hence expresses a preference to pure
(non-mixed) phases $\phi = \pm 1$. Also, the explicit penalty term for
steep gradients causes any interfaces to become diffused.

For $\psi$, a \emph{logarithmic} free energy is rather preferred as
this provides for certain analytical properties such as possessing an
\emph{isotherm} relation as outlined in Section~\ref{subsec:isotherm},
\begin{align}
  \label{eq:Fpsi}
  \Fpsi &= \Pii \left[ \psi\log \psi+(1-\psi)\log(1-\psi) \right].
\end{align}
The temperature-dependent constant $\Pii$ (cf.~\eqref{eq:newparams})
takes the role of a diffusion coefficient for $\psi$ and $\Fpsi$ is
thus the energy potential governing the entropy decrease of mixing the
surfactant with the bulk phase. We shall later find reasons to modify
\eqref{eq:Fpsi} since it, unlike \eqref{eq:Fphi} does not contain a
square gradient term. Note that a logarithmic free energy is to be
combined with a \emph{degenerate} mobility which vanishes at the
extreme points $\psi \in \{0,1\}$ \cite[Sect.~4.2]{Ch4CHE}. The usual
choice is simply
\begin{align}
  \label{eq:Mpsi}
  M_{\psi} = \psi(1-\psi)
\end{align}
which will also be used here. A point in favor of formulating the
Cahn-Hilliard equation using this kind of mobility is that, unlike the
constant mobility case, it can rigorously be shown to produce
solutions $0 \le \psi \le 1$ \cite[Sect.~4.2]{Ch4CHE}.

Surfactant molecules tend to move to a point and orient themselves in
such a way that they are at ease. Owing to their nature this generally
means at the interface between fluids, with the hydrophilic part
pointing towards a fluid of the system and the hydrophobic part
pointing away from this fluid. In the current model, $\Fone$ is the
surface energy potential accounting for this adsorption,
\begin{align}
  \label{eq:Fone}
  \Fone &= -\frac{\Cn^{2}}{4} \psi (\nabla \phi)^{2}.
\end{align}
This particular form is inspired by the corresponding term in the
sharp interface model of \cite{surfkinetics0,sharpsurf}. The square
gradient acts as a diffuse version of the sharp interface indicator
function and can rigorously be interpreted as a nascent Dirac delta
function as outlined in Section~\ref{subsubsec:Model2and3} below.

Finally, to conclude the specification of the free energy, the term
$\Fex$ takes the form
\begin{align}
  \label{eq:Fex}
  \Fex &= \frac{1}{4\Ex} \psi\phi^{2}
\end{align}
and penalizes free surfactant in the respective phases. Similarly to
$\Fpsi$, $\Fex$ is an enthalpic term measuring the cost of free
surfactant. However, due to the factor $\phi^{2}$, $\Fex$ is inactive
at an interface where $\phi \approx 0$. Actually, \eqref{eq:Fex} is
just the simplest form among many possibilities
(cf.~\cite{ternaryfluid, surfdynamics, diffemul}). To some extent,
\eqref{eq:Fone} and \eqref{eq:Fex} are complementary: $\Fone$ locally
attracts surfactant to an existing interface while $\Fex$ globally
counteracts the occurrence of free surfactant. Finally, and as we
shall see in Section~\ref{subsec:equilibrium}, \eqref{eq:Fex} can also
be regarded as defining the \emph{bulk solubility}.

Through variational derivatives of the free energy $F$ with respect to
the two variables $\phi$ and $\psi$ one now obtains the chemical
potentials
\begin{align}
  \label{eq:muphi}
  \mu_{\phi} &= \frac{\delta F}{\delta \phi} = -\phi+\phi^{3}-
  \frac{\Cn^{2}}{2} \Delta \phi+\frac{\Cn^{2}}{2} \psi \Delta \phi+
  \frac{\Cn^{2}}{2} \nabla \psi \cdot \nabla \phi+
  \frac{1}{2\Ex} \psi \phi, \\
  \label{eq:mupsi}
  \mu_{\psi} &=  \frac{\delta F}{\delta \psi} = 
  \Pii \log \frac{\psi}{1-\psi}-\frac{\Cn^{2}}{4}
  (\nabla \phi)^{2}
  +\frac{1}{4\Ex} \phi^{2}.
\end{align}

Note that, carrying out the differentiation in the right-hand side of
\eqref{eq:psiPDE} for the logarithmic term in \eqref{eq:mupsi}, we get
using \eqref{eq:Mpsi} that
\begin{align}
   \label{eq:mupsiR}
  \frac{1}{\Pe_{\psi}}
  \nabla \cdot M_{\psi} \nabla \Pii \left( \log \frac{\psi}{1-\psi}
  \right) &= \frac{\Pii}{\Pe_{\psi}} \Delta \psi,
\end{align}
provided that $0 < \psi < 1$.

This completes the specification of the \emph{gradient flow}
\eqref{eq:phiPDE}--\eqref{eq:psiPDE} which together with the
Navier-Stokes equations \eqref{eq:cont}--\eqref{eq:NS} and suitable
boundary conditions make up our basic multiphase-surfactant model. In
what follows we shall initially be concerned mainly with the
diffusion-controlled part \eqref{eq:phiPDE}--\eqref{eq:psiPDE} (with
$\fatu = 0$) and postpone computational experiments with the full
hydrodynamic set of equations until Section~\ref{sec:numexp2}. Testing
\eqref{eq:phiPDE} with $\mu_{\phi}$ and \eqref{eq:psiPDE} with
$\mu_{\psi}$, respectively, we find under natural boundary conditions
that
\begin{align}
  \label{eq:energydecay}
  \frac{d}{dt} \int_{\Omega} F \, dx &=
  -\int_{\Omega} \frac{1}{\Pe_{\phi}}\|\nabla \mu_{\phi}\|^{2}+
  \frac{M_{\psi}}{\Pe_{\psi}}\|\nabla \mu_{\psi}\|^{2} \, dx \le 0
\end{align}
for $M_{\psi} \ge 0$. This relation expresses the decay of the free
energy and is required for thermodynamical consistency. For the full
hydrodynamic set of equations we find similarly by testing
\eqref{eq:NS} with $\Ca\Cn\Reynold \fatu$ the general energy identity
\begin{align}
  \nonumber
  \frac{d}{dt} &\int_{\Omega}
  \frac{\rho\Ca\Cn\Reynold}{2}\|\fatu\|^{2}+F \, dx = \\
  \label{eq:energyid}
  &\phantom{=}
  -\int_{\Omega} \frac{1}{\Pe_{\phi}}\|\nabla \mu_{\phi}\|^{2}+
  \frac{M_{\psi}}{\Pe_{\psi}}\|\nabla \mu_{\psi}\|^{2}+
  \frac{\rho\nu\Ca\Cn}{2} \|\nabla \fatu+(\nabla \fatu)^{T}\|^{2}
  \, dx.
\end{align}


\section{Analysis}
\label{sec:analysis}

In this section we first show that the free-space PDE defined by
\eqref{eq:phiPDE}--\eqref{eq:psiPDE} with $\fatu = 0$ is ill-posed in
the sense of frozen coefficients in a neighborhood of a sufficiently
smooth steady-state solution. We then return to the baseline model
\eqref{eq:phiPDE}--\eqref{eq:NS} and modify it in order to provide for
well-posedness. Three different modifications are suggested, all with
differing motivations and properties.

\subsection{Ill-posedness}
\label{subsec:illposedness}

We thus consider the case $\fatu = 0$ in
\eqref{eq:phiPDE}--\eqref{eq:psiPDE} for a one-dimensional free-space
formulation. We make the quasi steady-state ansatz $\phi =
\phi_{\eq}+\delta u$ and $\psi = \psi_{\eq}+\delta v$, where
$\phi_{\eq}/\psi_{\eq}$ are equilibrium solutions as $t \to \infty$
and where $\delta$ is a small parameter. Assuming the existence of
equilibrium solutions is motivated by the fact that
\eqref{eq:phiPDE}--\eqref{eq:psiPDE} describe a \emph{gradient flow}
and naturally possesses (non-unique) equilibria
\cite[Sect.~4]{ODEstabrev}. ---Fundamentally, the lack of any
equilibrium would violate the second law of thermodynamics.

Assume for now that the equilibrium solutions are sufficiently regular
that a linearization in $\delta$ is possible (see the discussion after
\eqref{eq:lambda2} below). Upon linearizing and keeping only the
principal part of the operator we get after some work the linear
system of PDEs
\begin{align}
  \label{eq:ppart}
  \left[ \begin{array}{c}
      u_{t} \\
      v_{t}
    \end{array} \right] &= 
  \left[ \begin{array}{cc}
      -\frac{\Cn^{2}}{2}\frac{1-\psi_{\eeq}}{\Pe_{\phi}}D^{4} & 
      \frac{\Cn^{2}}{2}\frac{D(\phi_{\eeq})}{\Pe_{\phi}}D^{3} \\
      -\frac{\Cn^{2}}{2}\frac{\psi_{\eeq}(1-\psi_{\eeq})D(\phi_{\eeq})}
      {\Pe_{\psi}}D^{3}
      & \frac{\Pii}{\Pe_{\psi}} D^{2}
    \end{array} \right] \left[ \begin{array}{c}
      u \\
      v
    \end{array} \right],
\end{align}
where $D = d/dx$.  Without further ado we take the ``frozen
coefficient'' Fourier transform of \eqref{eq:ppart} by simply
exchanging $D \to -i\omega$,
\begin{align}
  \label{eq:ppart2}
  \left[ \begin{array}{c}
      \hat{u}_{t} \\
      \hat{v}_{t}
    \end{array} \right] &= 
  \left[ \begin{array}{cc}
      -\frac{\Cn^{2}}{2}\frac{1-\psi_{\eeq}}{\Pe_{\phi}}\omega^{4} & 
      \frac{\Cn^{2}}{2}\frac{D(\phi_{\eeq})}{\Pe_{\phi}}i\omega^{3} \\
      -\frac{\Cn^{2}}{2}\frac{\psi_{\eeq}(1-\psi_{\eeq})D(\phi_{\eeq})}
      {\Pe_{\psi}}i\omega^{3}
      & -\frac{\Pii}{\Pe_{\psi}} \omega^{2}
    \end{array} \right] \left[ \begin{array}{c}
      \hat{u} \\
      \hat{v}
    \end{array} \right].
\end{align}
Call the above matrix $A$. Firstly, since the characteristic equation
is real, complex eigenvalues must come in conjugate pairs. Secondly,
from a quick calculation we see that $\det A \propto \Pii-\Cn^{2}/2
\cdot \psi_{\eq}(D(\phi_{\eq}))^{2}$ and hence, when the determinant
is negative there is a positive real eigenvalue indicating some kind
of instability. Slightly more elaborate calculations are now required
to show that there is in fact an eigenvalue
\begin{align}
  \label{eq:lambda2}
  \lambda_{2} &= \frac{1}{\Pe_{\psi}} \left[ 
    \frac{\Cn^{2}}{2} \psi_{\eq}(D(\phi_{\eq}))^{2}-\Pii
  \right] \omega^{2}+\Ordo{1},
\end{align}
while $\lambda_{1}$ is asymptotically negative and scales as
$\Ordo{\omega^{4}}$. From \eqref{eq:lambda2} we see that $\lambda_{2}$
may well be positive for $\omega$ sufficiently large. If this is the
case the PDE violates the \emph{Petrovskii condition} which is
necessary for well-posedness \cite[Theorem~4.5.2]{blueBook}.

To turn these arguments into a strict proof unfortunately requires
non-trivial pointwise estimates of the equilibrium solutions. Besides
the fact that there exist problems for which the frozen coefficients
formulation is unstable while the variable coefficient problem is
stable, the linearization argument can also be questioned. Solutions
$\phi$ and $\psi$ of \eqref{eq:phiPDE}--\eqref{eq:psiPDE} could loose
derivatives as $t \to \infty$ such that $\phi_{\eq}$ and $\psi_{\eq}$
are sufficiently non-smooth that the linearization argument becomes
invalid. The only immediately available \textit{a priori} estimates
are that $\phi_{\eq} \in \Hone(\Omega)$ and that $\psi_{\eq} \in
\Ltwo(\Omega)$ (this follows directly from \eqref{eq:energydecay}). We
shall not delve more into these matters here, but only remark that a
linear stability analysis is well motivated from the viewpoint of
physical modeling and also, that all of our numerical experiments
agree very well with the analysis performed above
(cf.~Figure~\ref{fig:isotherm0}).

Using results to be developed in Section~\ref{sec:numexp1} we now
derive some simplified conditions for the instability to occur. Let us
consider it inherent in the model that both $\psi_{\eq}$ and
$(D(\phi_{\eq}))^{2}$ are largest at the interface (taken to be at $x
= 0$). Then $\psi_{\eq} \le \psi_{\eq}(0) = \psi_{0}$, the surfactant
loading at the interface and, provided that the interface sharpness is
independent of the surfactant loading, $(D(\phi_{\eq}))^{2} \le
(D(\phi_{\eq}(0)))^{2} \sim 1/\Cn^{2}$ (see
Section~\ref{subsec:equilibrium}). Using the isotherm relation
$\psi_{0} \sim \psi_{b}/(\psi_{c}+\psi_{b})$ (again, see
Section~\ref{subsec:isotherm} below), an \textit{a priori} approximate
sufficient condition for \emph{in}stability is therefore
\begin{align}
  \label{eq:Pilim}
  \Pii &< \frac{\psi_{b}/2}{\psi_{c}+\psi_{b}}, \\
  \intertext{expressed in terms of the \emph{bulk value} $\psi_{b} =
    \psi_{\eq}(\pm \infty)$ and the \emph{Langmuir adsorption
      constant} $\psi_{c}$ discussed thoroughly below. Using the
    relation \eqref{eq:psic} for the latter this is equivalent to}
  \label{eq:psiblim}  
  \psi_{b} &> \frac{2\Pii}{1-2\Pii} \times \psi_{c} = 
  \frac{2\Pii}{1-2\Pii} \times 
  \exp \left( -\frac{1+1/\Ex}{4\Pii} \right)
\end{align}
provided that $\Pii < 1/2$. That is, sufficiently large values of
$\psi_{b}$ lead to an instability near the equilibrium solution (which
therefore does not exist).

This situation is somewhat remindful of results for a close relative
of the Cahn-Hilliard equation, namely the \emph{Thin film
  equation}. This model can also be derived from a Ginzburg-Landau
energy which, as in \eqref{eq:energydecay}, by construction decays
with time. Through a clever constructive argument, however, one can
show that for a sufficiently large total mass, there is a smooth
initial profile for which finite-time blowup in the form of indefinite
focusing occurs (see \cite{unstable_thinfilm} and the references
therein). Unfortunately, it seems difficult to transfer these types of
arguments into the current setting.

\begin{example}
  In several numerical experiments in \cite{diffsurf} using the
  lattice Boltzmann method, $\Pii \in [0.08,0.20]$ and $\Ex \sim 1$.
  To be concrete and still use actual (published) values, we take
  $\Pii = 0.1227$ with $\Ex = 1$. Then \eqref{eq:psiblim} becomes
  $\psi_{b} > 5.526 \times 10^{-3}$ and under this condition we expect
  the model to be ill-posed. Since many of the experiments in
  \cite{diffsurf} (see for example Fig.~1 and 2) are made in this
  regime, these results are misleading at the very least. One cannot
  help speculating that the lattice Boltzmann method used enjoy
  ``hidden'' higher order derivative terms that help increasing the
  stability. In Figure~\ref{fig:unstable}, p.~\pageref{fig:unstable}
  these predictions are illustrated using the numerical method
  developed in Appendix~\ref{app:num1}. We note that similar
  experiments are made anew in \cite{diffsurf2} under only slightly
  altered conditions. Here, the parameters in the simulations appear
  to have been kept within the limits for well-posedness. Using this
  model substantially limits the physical parameter space that can be
  explored.
\end{example}

Incidentally, \eqref{eq:psiblim} disqualifies the label ``numerical
stabilization term'' for the bulk-energy $\Fex$ in \eqref{eq:Fex} as
used in \cite[p.~4]{diffsurf} and again in \cite[p.~9167]{diffsurf2}.
Increasing their variable $W$ in \eqref{eq:newparams} (that is,
``increasing'' the numerical stabilization) is the same thing as
\emph{decreasing} $\Ex$, thereby driving more surfactant to the
interface. Assuming that the other parameters $\Cn$ and $\Pii$ are
kept fixed, according to \eqref{eq:psiblim} this means that
\emph{less} amount of surfactant implies \emph{mathematical} (as
opposed to \emph{numerical}) instability.

Finally, and from the perspective of physical modeling, one could in
principle accept the ill-posedness as an artifact of the fact that we
are dealing with an artificially diffuse interface. However, in such a
case it is natural to demand at the very least that the limiting
problem $\Cn \to 0$ is stable. For the eigenvalue $\lambda_{2}$ in
\eqref{eq:lambda2} to be negative as $\Cn \to 0$ one clearly has to
give up the natural requirement that the interface sharpness is
independent of the surfactant loading. Indeed, our approximate bounds
\eqref{eq:Pilim}--\eqref{eq:psiblim} are fully independent of $\Cn$ as
a reflection of this assumption.

\subsection{Three alternatives}
\label{subsec:alt}

In this section we seek viable alternatives to the model just
analyzed. Our first modified formulation involves using the full form
of the logarithmic free energy $\Fpsi$ in \eqref{eq:Fpsi}. This form
contains a square gradient term and leads to an additional fourth
order dissipative term in the PDE for the surfactant concentration
$\psi$. The downside is that the resulting model has more complicated
mathematical properties and we therefore look for simpler but still
reasonable alternatives. Two such suggestions are obtained by using a
derivative-free diffuse interface Dirac delta function when defining
$\Fone$ in \eqref{eq:Fone}.

\subsubsection{Completing the logarithmic free energy}
\label{subsubsec:Model1}

Rather than \eqref{eq:Fpsi}, the complete logarithmic free energy is
actually given by (\cite[Sect.~4.2, Eqs.~(20), (21)]{Ch4CHE}; see also
\cite{multiphasefield})
\begin{align}
  \label{eq:Fpsi1}
  \Fpsi &= \Pii \left[ \psi\log \psi+(1-\psi)\log(1-\psi)
  \right]+\frac{\sigma}{4}\psi(1-\psi)+\frac{\Cn^{2}}{4} (\nabla \psi)^{2},
\end{align}
with $\sigma$ a new parameter. Accordingly, our \emph{``Model 1''} is
obtained by exchanging the definition of $\Fpsi$ in \eqref{eq:Fpsi}
with \eqref{eq:Fpsi1}. While the chemical potential $\mu_{\phi}$ in
\eqref{eq:muphi} is not affected by this addition to the free energy,
$\mu_{\psi}$ in \eqref{eq:mupsi} is to be replaced with
\begin{align}
  \label{eq:mupsi1}
  \mu_{\psi} &= \Pii \log \frac{\psi}{1-\psi}-\frac{\Cn^{2}}{4}
  (\nabla \phi)^{2}-\frac{\Cn^{2}}{2} \Delta \psi-\frac{\sigma}{2}\psi
  +\frac{1}{4\Ex} \phi^{2}.
\end{align}
Evidently, one more boundary condition than before is needed since now
the equation for $\psi$ is fourth order.

A possible explanation to the fact that the square gradient term is
missing in the original presentation of the model in \cite{diffsurf}
is that it is also omitted in the earlier references
\cite{surfkinetics0,sharpsurf,surfdynamics}. In
\cite{surfkinetics0,sharpsurf} we note that the term has no meaning
since the interface model there is sharp. Interestingly, the square
gradient term \emph{is} in fact included in the surfactant models
found in \cite{ternaryfluid,dropletsurf,diffemul}, but then for the
constant mobility case. As remarked in \cite{diffsurf} (and as will be
evident in Section~\ref{subsec:isotherm} below), such models do not
possess a realistic adsorption isotherm.

This also seems to be the place to mention how the two order
parameters $(\phi,\psi)$ are to be interpreted. Rather than a true
multicomponent system as treated in \cite{multiphasefield}, we are
dealing with small, often \emph{very} small, concentrations of
surfactant. A kind of two-step approximation is therefore employed
where a Cahn-Hilliard equation with degenerate mobility governs the
surfactant ($\psi$) and ``non-surfactant'' ($1-\psi$) phases. In turn,
the non-surfactant phase (e.g.~oil and water) is treated by a
classical constant mobility Cahn-Hilliard equation in the phase-field
variable $\phi$. The three concentration components of the flow should
therefore rightly be understood as
\begin{align}
  &[\psi,(1-\psi)(1+\phi)/2,(1-\psi)(1-\phi)/2]
  \intertext{and are thus approximated with}
  \label{eq:phasesapprox}
  &[\psi,(1+\phi)/2,(1-\phi)/2]
\end{align}
so that the volume is only preserved up to $\Ordo{\psi}$.

We conclude by offering a few comments on the \emph{lateral
  interaction term} $-\sigma/4 \cdot \psi^{2}$ in \eqref{eq:Fpsi1}.
This term is included also in \cite{surfkinetics0,sharpsurf,diffsurf2}
where it is shown to lead to the \emph{Frumkin isotherm} (see
Section~\ref{subsec:isotherm}). Interestingly, the same term appears
also in \cite{ternaryfluid,surfdynamics,dropletsurf,diffemul}, but
here with the opposite sign. Clearly, a positive value of $\sigma$
favors \emph{clustering} of surfactant since the term can be
interpreted as a free energy decrease for \emph{pairs} of surfactant
molecules. By differentiating twice the non-gradient part of
\eqref{eq:Fpsi1} with respect to $\psi$ we further find that
\begin{align}
  \label{eq:sigmaconvex}
  \frac{d^{2}}{d\psi^{2}} \left[ \Pii \left[ \psi\log \psi+(1-\psi)\log(1-\psi)
    \right]+\frac{\sigma}{4}\psi(1-\psi) \right] &= 
  \frac{\Pii}{\psi(1-\psi)}-\frac{\sigma}{2}.
\end{align}
Hence, for $0 < \psi < 1$, this part of the free energy remains convex
in $\psi$ provided that we choose $\sigma \le 8\Pii$.

\subsubsection{Gradient-free Dirac delta functions}
\label{subsubsec:Model2and3}

We have already commented that the term $-\Cn^{2}/4 \cdot
\psi(\nabla\phi)^{2}$ in \eqref{eq:Fone} results from using the square
gradient as a diffuse version of the sharp interface indicator
function. Indeed, for the planar interface equilibrium solution
$\phi(x) = \tanh (x/\Cn)$, we have that $(\nabla\phi)^{2} = \Cn^{-2}
\sech^{4}(x/\Cn)$. Over the real line, the quartic hyperbolic secant
is a nascent Delta function in the sense that $3/(4\Cn)
\sech^{4}(x/\Cn) \to \delta(x)$ (convergence in distribution in the
sharp interface limit $\Cn \to 0$).

Noting that the function $1-\tanh^{2}(x/\Cn) = \sech^{2}(x/\Cn)$ also
defines a nascent Delta function by virtue of the limit $1/(2\Cn)
\sech^{2}(x/\Cn) \to \delta(x)$, a tempting replacement for
$(\nabla\phi)^{2}$ in \eqref{eq:Fone} is $2/(3\Cn^{2}) \cdot
(1-\phi^{2})$ after appropriate scaling. However, as we shall see,
this choice has a different behavior near the origin and it shall
later be convenient to use the slightly altered scaling $1/\Cn^{2}
\cdot (1-\phi^{2})$ (see Section~\ref{subsec:isotherm}). We thus
define our \emph{``Model 2''} by replacing \eqref{eq:Fone} with
\begin{align}
  \label{eq:Fone2}
  \Fone &= -\frac{1}{4} \psi (1-\phi^{2}).
\end{align}
Taking variational derivatives we obtain the new chemical potentials
\begin{align}
  \label{eq:muphi2}
  \mu_{\phi} &= -\phi+\phi^{3}-
  \frac{\Cn^{2}}{2} \Delta \phi+\frac{1}{2}\psi\phi+
  \frac{1}{2\Ex} \psi \phi, \\
  \label{eq:mupsi2}
  \mu_{\psi} &= \Pii \log \frac{\psi}{1-\psi}+\frac{1}{4}\phi^{2}
  +\frac{1}{4\Ex} \phi^{2},
\end{align}
where it is clear that Model 2 can be implemented by simply skipping
the term $\Fone$ altogether and substitute $\Ex \to 1/(1/\Ex+1)$ in
\eqref{eq:Fex}.

Although the square hyperbolic secant is a well-known nascent Delta
function (it is the derivative of the Fermi-Dirac function), in the
present context one can actually continue to use the quartic
hyperbolic secant, without introducing any derivatives. A simple
replacement for the square gradient achieving just this is $1/\Cn^{2}
\cdot (1-\phi^{2})^{2}$. Accordingly, our \emph{``Model 3''} uses, in
place of \eqref{eq:Fone},
\begin{align}
  \label{eq:Fone3}
  \Fone &= -\frac{1}{4} \psi (1-\phi^{2})^{2},
\intertext{with the associated chemical potentials}
  \label{eq:muphi3}
  \mu_{\phi} &= -\phi+\phi^{3}-
  \frac{\Cn^{2}}{2} \Delta \phi+(1-\phi^{2})\psi\phi+
  \frac{1}{2\Ex} \psi \phi, \\
  \label{eq:mupsi3}
  \mu_{\psi} &= \Pii \log \frac{\psi}{1-\psi}-
  \frac{(1-\phi^{2})^{2}}{4}+\frac{1}{4\Ex} \phi^{2}.
\end{align}
Evidently, Model 3 is stiffer than Model 2 and it will also be shown
to produce a sharper equilibrium profile. Figure~\ref{fig:unstable}
and \ref{fig:stable} below show representative sample simulations of
all four models thus far considered. These results were obtained with
the numerical method discussed in Appendix~\ref{app:num1}.

\subsection{Sample simulation: ill-posedness}

\begin{figure}
  \ifpdf
  \includegraphics[width = 149pt]{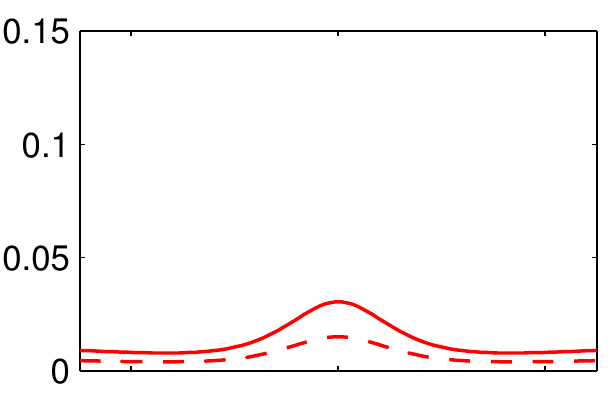}
  \includegraphics[width = 149pt]{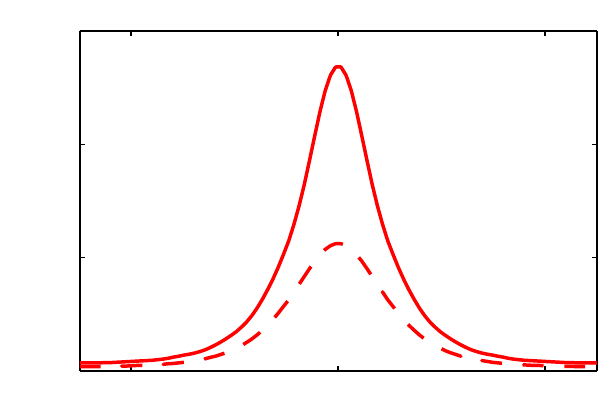}
  \includegraphics[width = 149pt]{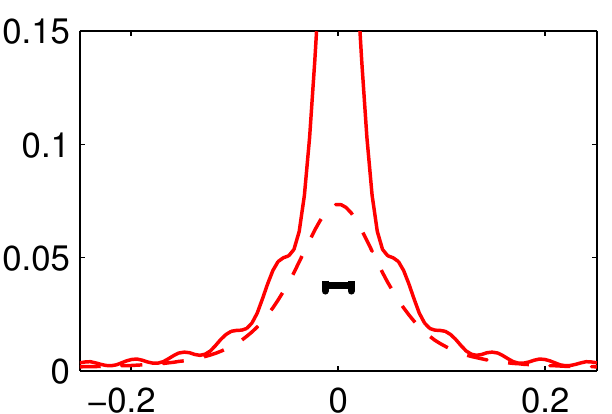}
  \includegraphics[width = 149pt]{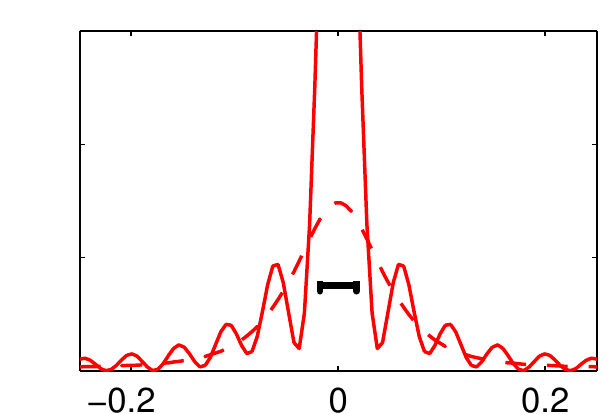}
  \includegraphics[width = 305pt]{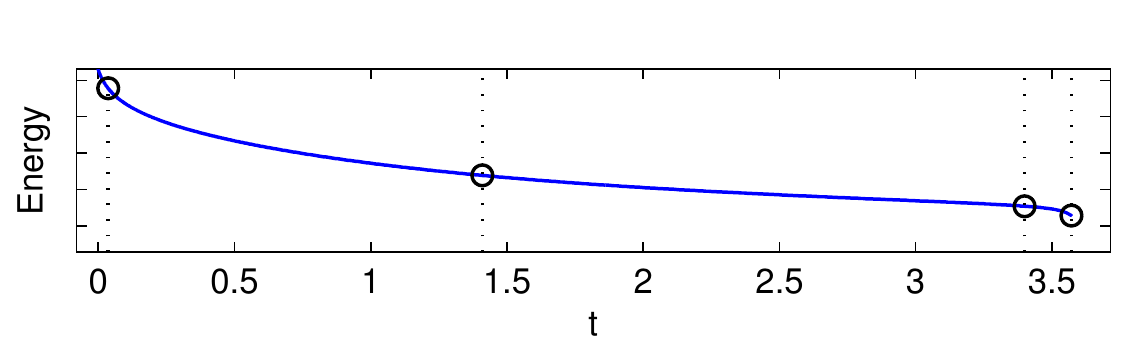}
  \else
  \includegraphics[width = 149pt]{fig/unstable1.eps}
  \includegraphics[width = 149pt]{fig/unstable2.eps}
  \includegraphics[width = 149pt]{fig/unstable3.eps}
  \includegraphics[width = 149pt]{fig/unstable4.eps}
  \includegraphics[width = 305pt]{fig/unstable_energy.eps}
  \fi
  \caption{Time snapshots of the surfactant concentration for an
    unstable \textit{(solid)} and a stable \textit{(dashed)} case. All
    parameters are the same in the two cases ($\Cn = 1/6$, $\Ex = 1$,
    and $\Pii = 0.1227$), but the amount of surfactant differs by a
    factor of two. The initial data for $\psi(x)$ is a flat profile
    with height 0.012 and 0.006, respectively, and $\phi(x)$ is
    initially and in both cases set to $\tanh(x/\Cn)$. In the top two
    graphs the instability has not yet developed. In the two graphs
    below the $\Ordo{\omega^{2}}$-part of \eqref{eq:lambda2} is
    positive inside the indicated small region near the origin, and
    the instability immediately becomes manifest. In the bottom graph
    the free energy \eqref{eq:F}--\eqref{eq:Fex} is plotted as a
    function of time $t$ and the times for the four snapshots are
    indicated by circles. See text for further comments.}
  \label{fig:unstable}
\end{figure}

In Figure~\ref{fig:unstable} we compare an unstable and a stable case
of the original surfactant phase-field model (referred to as ``Model
0'' from now on). The initial data, in the form of a uniform
surfactant profile, was chosen in accordance with the stability
criterion \eqref{eq:psiblim} and in such a way that both cases are in
close proximity to the boundary of the region of
well-posedness. During the simulation time displayed here it holds for
the numerical solution $\psi$ that $0 < \psi < 1$ and that the
associated Ginzburg-Landau \eqref{eq:F} energy is decreasing. Shortly
after the displayed simulation time the numerical solution becomes
negative such that the free energy formally becomes multivalued.

In Section~\ref{subsec:isotherm} we more carefully evaluate the
sharpness of the condition \eqref{eq:psiblim} after developing some
more concepts (see Figure~\ref{fig:isotherm0}).

\begin{figure}
  \ifpdf
  \includegraphics[width = 149pt]{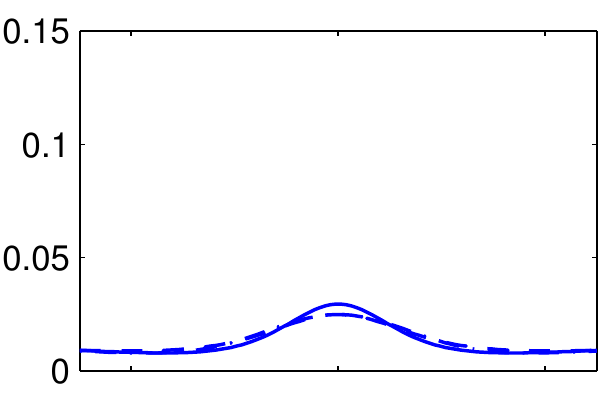}
  \includegraphics[width = 149pt]{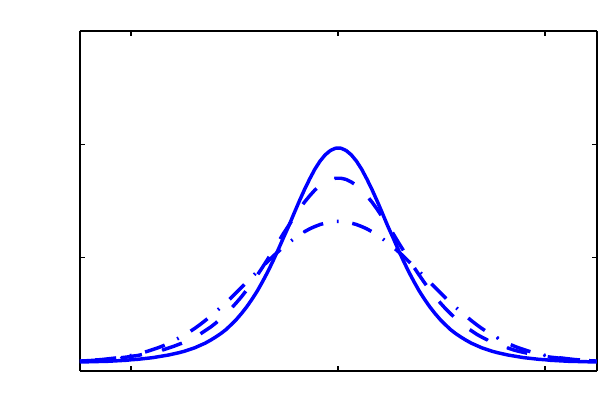}
  \includegraphics[width = 149pt]{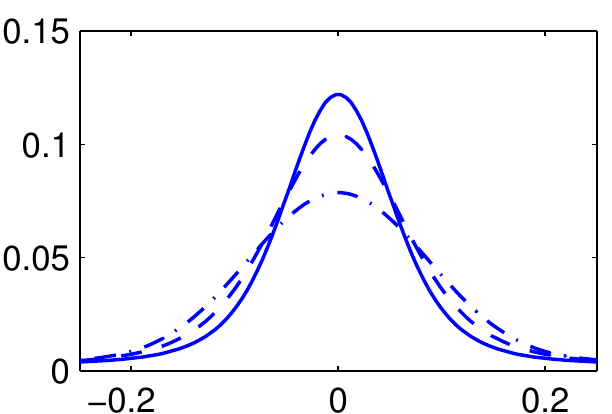}
  \includegraphics[width = 149pt]{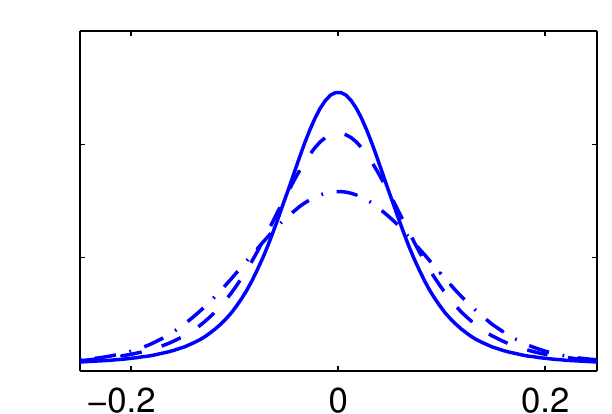}
  \else
  \includegraphics[width = 149pt]{fig/stable1.eps}
  \includegraphics[width = 149pt]{fig/stable2.eps}
  \includegraphics[width = 149pt]{fig/stable3.eps}
  \includegraphics[width = 149pt]{fig/stable4.eps}
  \fi
  \caption{Snapshots and parameters as in the unstable case of
    Figure~\ref{fig:unstable}, but using Model 1--3
    instead. \textit{Dash-dot:} Model 1 with $\sigma = 8\Pii$,
    \textit{dashed:} Model 2, \textit{solid:} Model 3.}
  \label{fig:stable}
\end{figure}

Examples of all alternative models proposed in
Section~\ref{subsec:alt} are displayed in Figure~\ref{fig:stable} were
we simulate the ill-posed case from Figure~\ref{fig:unstable}
anew. Unlike Model 0, for this choice of parameters, all three new
models are perfectly stable.

\section{Analysis and numerical experiments in 1D}
\label{sec:numexp1}

In order to evaluate and compare the proposed models we now proceed to
derive some analytical estimates, most of which we test through
numerical simulations. We discuss equilibrium solutions, adsorption
isotherm relations, and we also evaluate experimentally the
diffusion-controlled adsorption dynamics at the interface.

\subsection{Planar equilibrium solution at constant surfactant
  concentration}
\label{subsec:equilibrium}

If $\psi$ is held constant at a bulk value $\psi_{b}$, steady-state of
$\phi$ implies for both our starting model \eqref{eq:muphi} and for
Model 1 in Section~\ref{subsubsec:Model1} that
\begin{align}
  \mu_{\phi} &= -\phi+\phi^{3}-
  \frac{\Cn^{2}}{2} \Delta \phi+\frac{\Cn^{2}}{2} \psi_{b} \Delta \phi+
  \frac{1}{2\Ex} \psi_{b} \phi = 0.
\end{align}
Generally, steady-state requires $\mu_{\phi} = \mbox{constant}$, but
it is not difficult to see that a solution which is odd about the
origin (say, $\phi(\pm \infty) = \pm (1-\varepsilon)$) must in fact
have $\mu_{\phi} = 0$. Rewriting we get
\begin{align}
  \nonumber
  \mu_{\phi} &= -\underbrace{\left(1-\frac{1}{2\Ex} \psi_{b}\right)}_
  {=:\phi_{b}^{2}} \phi+\phi^{3}-
  (1-\psi_{b})\frac{\Cn^{2}}{2} \Delta \phi \\
  &= \phi_{b}^{3} \Bigl[ -\phi_{b}^{-1} \phi+(\phi_{b}^{-1} \phi)^{3}-
  \underbrace{\frac{1-\psi_{b}}{\phi_{b}^{2}}\frac{\Cn^{2}}{2}}_
  {=: \Cn_{\varphi}^{2}/2} \Delta \phi_{b}^{-1}\phi \Bigr].
\end{align}
This is just a new Cahn-Hilliard steady-state equation in the scaled
variable $\varphi := \phi_{b}^{-1}\phi$ with Cahn-number
$\Cn_{\varphi} = \Cn \sqrt{1-\psi_{b}}/\phi_{b}$. From the classical
solution $\varphi = \tanh (x/\Cn_{\varphi})$ we therefore get
\begin{align}
  \label{eq:phiinf0and1}
  \phi(x) &= \phi_{b} \tanh 
  \left( \phi_{b} \frac{x}{\Cn \sqrt{1-\psi_{b}}} \right), \\
\intertext{with}
  \label{eq:phib0and1}
  \phi_{b}^{2} &= 1-\frac{1}{2\Ex} \psi_{b}.
\end{align}
At the interface we have for this solution that $\Cn d/dx \, \phi(0) =
1-(1/\Ex-1)\psi_{b}/2+\Ordo{\psi_{b}^{2}}$ implying that the sharpness
of the interface is independent of the surfactant loading up to
$\Ordo{\psi_{b}^{2}}$ when $\Ex = 1$. The bulk behavior as $x \to \pm
\infty$ is $\pm \phi_{b}$ which means that the coefficient $\Ex$
controls the bulk solubility of the two phases relative to each
other. Namely, with $\phi_{b} = 1-\psi_{b}/(4\Ex)+\Ordo{\psi_{b}^{2}}$
(from \eqref{eq:phib0and1}), we see that \eqref{eq:phasesapprox}
becomes
\begin{align}
  &[\psi_{b},\psi_{b}/(8\Ex),1-\psi_{b}/(8\Ex)]
\end{align}
to within $\Ordo{\psi_{b}^{2}}$ in the two phases. Hence for a system
described by these thermodynamic potentials, the presence of the
surfactant in the bulk will allow for a small but finite solubility of
the one fluid in the other.

Proceeding in a similar fashion with Model 2 using \eqref{eq:muphi2}
we obtain the equilibrium solution
\begin{align}
  \phi(x) &= \phi_{b} \tanh \left( \phi_{b} \frac{x}{\Cn} \right), \\
  \phi_{b}^{2} &= 1-\left(1+\frac{1}{\Ex} \right) \frac{\psi_{b}}{2},
\end{align}
and where $\Cn d/dx \, \phi(0) = \phi_{b}^{2} = 1-\Ordo{\psi_{b}}$ and
is $< 1$ for any positive value of $\Ex$.

Finally, for Model 3 using \eqref{eq:muphi3} we get after similar
manipulations
\begin{align}
  \label{eq:phiinf3}
  \phi(x) &= \phi_{b} \tanh \left( \phi_{b} 
    \frac{x}{\Cn \sqrt{1-\psi_{b}}} \right), \\
\label{eq:phib3}
  \phi_{b}^{2} &= \frac{1-\left(1+\frac{1}{2\Ex} \right)
    \psi_{b}}{1-\psi_{b}}.
\end{align}
Not only are \eqref{eq:phiinf3} and \eqref{eq:phiinf0and1} identical
in form, but \eqref{eq:phib3} and \eqref{eq:phib0and1} in fact also
agree up to order $\Ordo{\psi_{b}^{2}}$. Consequently we have again
that $\Cn d/dx \, \phi(0) = 1-(1/\Ex-1)\psi_{b}/2+\Ordo{\psi_{b}^{2}}$
so that for $\Ex = 1$, the sharpness of the interface is independent
of the surfactant loading up to $\Ordo{\psi_{b}^{2}}$.

\begin{conclusion}
  For all four models, the parameter $\Ex$ controls the \emph{bulk
    solubility}. Model 0, 1, and 3 all agree closely in terms of the
  bulk value $\phi_{b}$ and also in the sharpness of the phase-field
  interface $\phi = 0$. By comparison, Model 2 has a more diffuse
  interface.
\end{conclusion}

\subsection{Adsorption isotherm}
\label{subsec:isotherm}

We now turn our attention to equilibrium profiles for $\psi$,
following the line of reasoning in \cite{diffsurf} closely. Since at
steady-state one must have that the chemical potential is constant,
the basic approach is to solve the equation $\mu_{\psi(x)} =
\mu_{\psi_{b}}$ for $\psi$, where $\psi_{b}$ as before is the bulk
concentration $\psi(\infty)$.  We write the chemical potentials for all
four models in the forms
\begin{align}
  \mu_{\psi_{b}} &= \Pii \log \frac{\psi_{b}}{1-\psi_{b}}+B_{b}
  +\frac{1}{4\Ex} \phi_{b}^{2}, \\
  \mu_{\psi(x)} &= \Pii \log \frac{\psi}{1-\psi}+B
  +\frac{1}{4\Ex} \phi^{2}. \\
  \intertext{Subtracting and introducing the intermediate variable
    $\psi_{c}(x)$ we get the relation}
  \Pii \log \psi_{c}(x) &= B-B_{b}-
  \frac{1}{4\Ex} \left(\phi_{b}^{2}-\phi^{2} \right) \\
  \intertext{in terms of which the steady-state profile is given by}
  \label{eq:psigen}
  \psi(x) &= \frac{\psi_{b}}{\psi_{b}+\psi_{c}(x)(1-\psi_{b})} =
  \frac{\psi_{b}}{\psi_{b}+\psi_{c}(x)}+\Ordo{\psi_{b}}.
\end{align}
The difference $B-B_{b}$ is given by, respectively,
\begin{align}
  \label{eq:dB}
  B-B_{b} &= \left\{
  \begin{array}{lr}
    -\frac{\Cn^{2}}{4}(\nabla\phi)^{2} &\mbox{(Model 0)} \\
    -\frac{\Cn^{2}}{4}(\nabla\phi)^{2}-
    \frac{\Cn^{2}}{2} \Delta\psi-\frac{\sigma}{2}(\psi-\psi_{b})
    &\mbox{(Model 1)} \\
    -\frac{1}{4}(\phi_{b}^{2}-\phi^{2}) &\mbox{(Model 2)} \\
    -\frac{1}{4}\left[(\phi_{b}^{2}-\phi^{2})
      (2-\phi_{b}^{2}-\phi^{2})\right] &\mbox{(Model 3)}
  \end{array} \right..
\end{align}
Clearly, for Model 1, this line of reasoning needs to be augmented
with additional assumptions or estimates in order not to be
circular. For the other cases, with a fixed phase-field profile
$\phi(x) \approx \phi_{b} \tanh (x/\Cn)$ (as determined in
Section~\ref{subsec:equilibrium}), \eqref{eq:psigen} yields a quite
decent approximation as we shall see.

Specializing $x$ in \eqref{eq:psigen} to the origin we get
\begin{align}
  \label{eq:isotherm}
  \psi_{0} &= \frac{\psi_{b}}{\psi_{b}+\psi_{c}}+\Ordo{\psi_{b}}, \\
  \label{eq:psic}
  \Pii \log \psi_{c} &= -\frac{1}{4} \left(1+\frac{1}{\Ex}\right)+
  \Ordo{\psi_{b}}.
\end{align}
To arrive at \eqref{eq:psic}, for Model 0 we have to assume that
$\Cn^{2}(\nabla\phi(0))^{2} = 1+\Ordo{\psi_{b}}$
(cf.~Section~\ref{subsec:equilibrium}), while for Model 2 and 3 we
only need to use the fact that $\phi(0) = 0$.

Eq.~\eqref{eq:isotherm} is the \emph{Langmuir isotherm} and $\psi_{c}$
as defined by \eqref{eq:psic} is the \emph{Langmuir (equilibrium)
  adsorption constant}. The fact that all models except Model 1
possesses the same Langmuir isotherm makes them comparable and also
explains our special choice of scaling when constructing Model 2 in
Section~\ref{subsubsec:Model2and3}.

For completeness, we note that \eqref{eq:psigen} for Model 1 becomes
as $x \to 0$,
\begin{align}
  \label{eq:isothermX}
  \psi_{0} &= \frac{\psi_{b}}{\psi_{b}+\psi_{c}
    \exp(-\alpha\psi_{0})R}+
  \Ordo{\psi_{b}},
\end{align}
with $\psi_{c}$ still defined by \eqref{eq:psic}, and where $\alpha =
\sigma/(2\Pii)$, $R = \exp(-\Cn^{2}\Delta\psi_{0}/(2\Pii))$. With $R =
1$, \eqref{eq:isothermX} is the \emph{Frumkin isotherm} which
essentially is an effect of the presence of the lateral interaction
term $-\sigma/4 \cdot \psi^{2}$ in \eqref{eq:Fpsi1}. Unfortunately,
there is no evident relation between $\psi_{0}$ and $\Delta \psi_{0}$
which can be used to close this line of reasoning. Therefore, for
Model 1 the adsorption isotherm cannot be explicitly determined.

\subsubsection{Numerical isotherm}
\label{subsubsec:num2}

Using our one-dimensional spectral Galerkin code as outlined in
Appendix~\ref{app:num1} we have performed several numerical
experiments with the adsorption isotherm relation. For all experiments
in this section we used the same set of parameters: $\Cn = 1/6$, $\Ex
= 1$, $\psi_{c} \in \{0.0020,0.0056,$ $0.016,0.035,0.075\}$ (with
$\Pii$ determined from the relation \eqref{eq:psic}), and with
$\psi_{b}$ sampled in the interval $[10^{-3},10^{-1}]$. All parameters
have been chosen to agree with those in \cite[Fig.~2]{diffsurf}. The
simulations were started with the profile $\phi(x) = \tanh(x/\Cn)$
with $\psi(x)$ defined by \eqref{eq:psigen} and, for convenience,
$\phi_{b} \equiv 1$. Finally, for Model 1 we used the value $\sigma =
8\Pii$, obtained from the requirement that the non-gradient part of
$\Fpsi$ be convex (see \eqref{eq:sigmaconvex}) and otherwise simply by
trial and error to approximately match the visual appearances of the
profiles for Model 2 and 3 (cf.~Figure~\ref{fig:samples}).

\begin{figure}
  \ifpdf
  \includegraphics[width = 301pt]{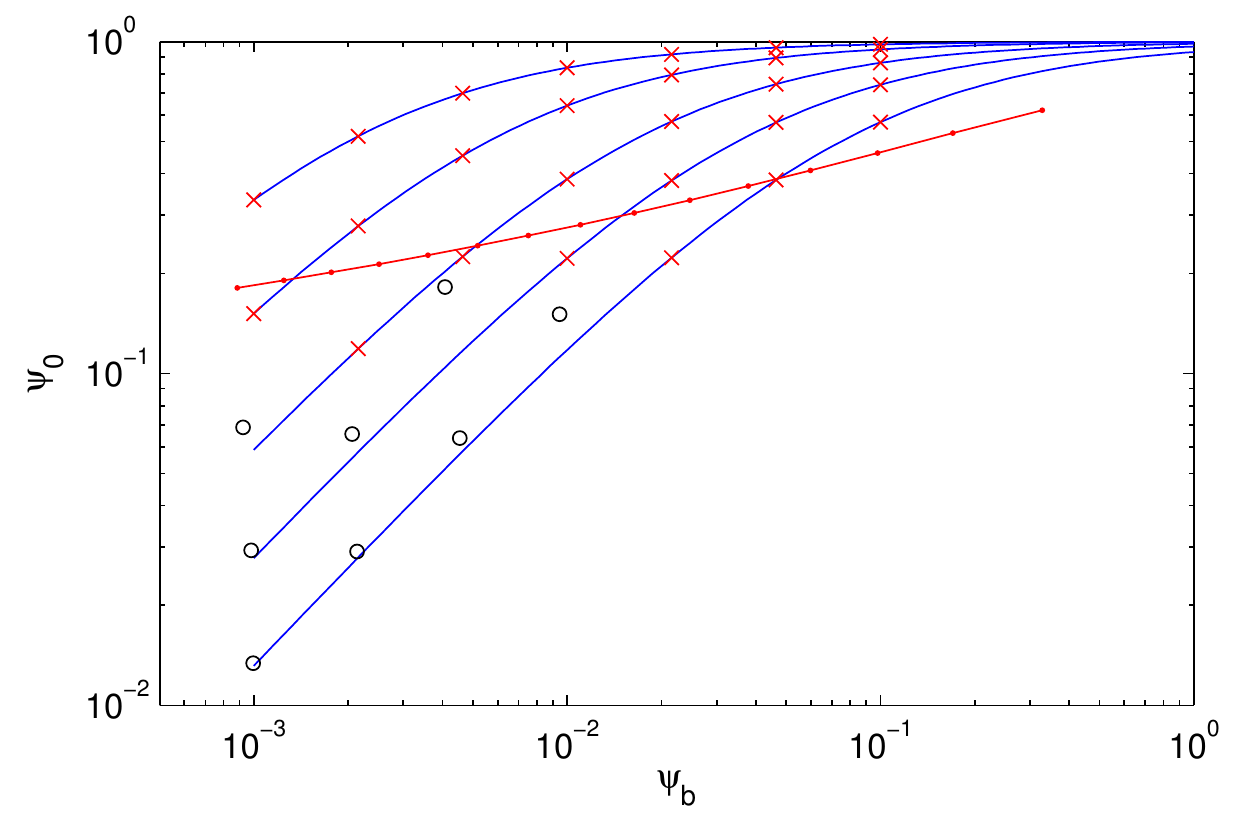}
  \else
  \includegraphics[width = 301pt]{fig/isotherm0.eps}
  \fi
  \caption{Attempt to verify the Langmuir isotherm \eqref{eq:isotherm}
    through numerical simulations of Model 0. This figure is intended
    to be an exact reproduction of the results in
    \cite[Fig.~2]{diffsurf}. However, here we clearly see that the
    model is ill-posed as predicted by the theory in
    Section~\ref{subsec:illposedness}. \textit{Solid:} the isotherm
    \eqref{eq:isotherm} for different values of $\psi_{c}$,
    \textit{circles:} numerical values. The \textit{dotted} line is
    the sufficient condition for ill-posedness \eqref{eq:psiblim} and
    \textit{crosses} are used to denote (missing) unphysical solutions
    (e.g.~large negative values).}
  \label{fig:isotherm0}
\end{figure}

In Figure~\ref{fig:isotherm0} we numerically test the sharpness of our
analysis of Model 0 in Section~\ref{subsec:illposedness}. Clearly, the
model is ill-posed for most of the parameters tested here, rendering
the model very questionable. From the figure, it is also seen that the
sufficient condition for ill-posedness \eqref{eq:psiblim} is quite
sharp.

In Figure~\ref{fig:isotherm0} a numerical isotherm is obtained for
Model 1. The Langmuir isotherm is not a bad model for small values of
$\psi_{c}$, but the measured adsorption breaks off for larger values
and also for higher surfactant concentrations. A few cases were found
to be numerically unstable by producing an extremely sharp profile for
$\phi$. We have not been able to analyze this and do not know at
present if this is an actual property of the model or a numerical
artifact of some kind.

\begin{figure}
  \ifpdf
  \includegraphics[width = 301pt]{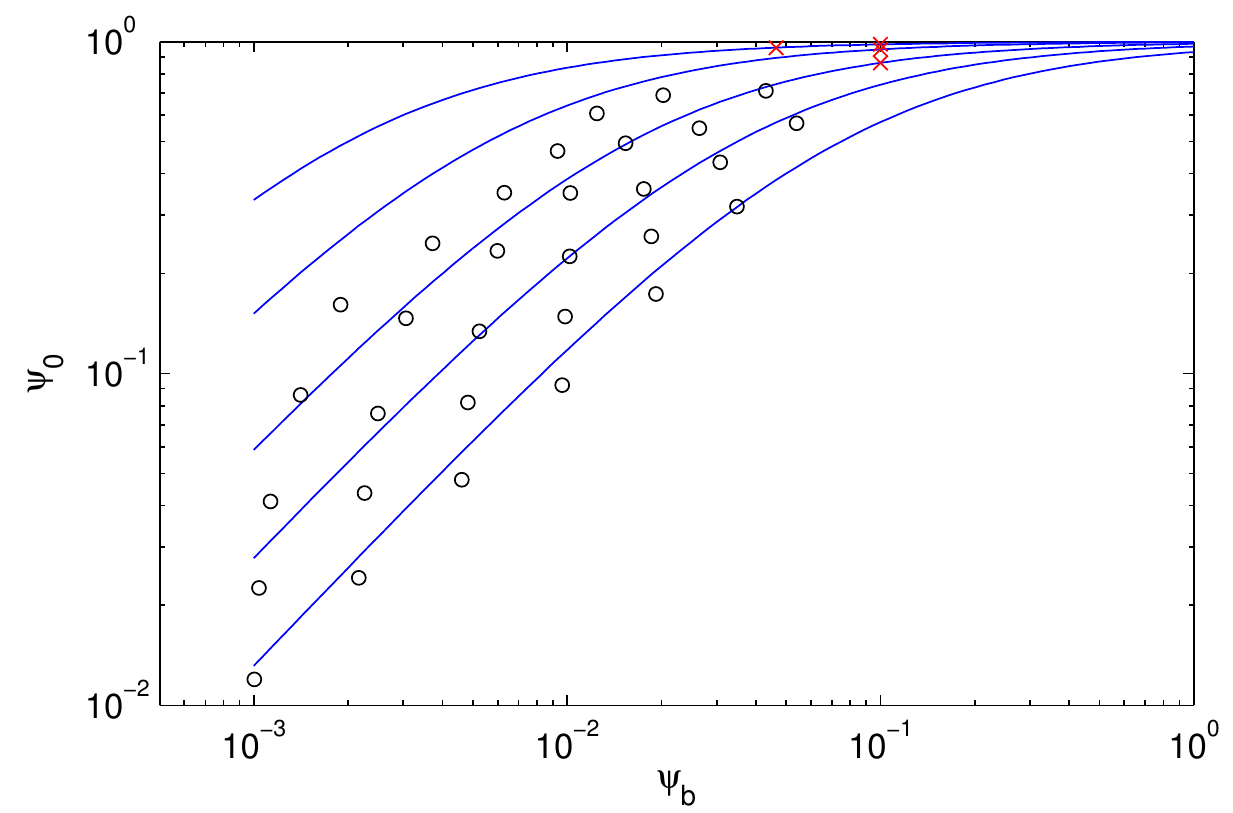}
  \else
  \includegraphics[width = 301pt]{fig/isotherm1.eps}
  \fi
  \caption{Results from simulations with parameters as in
    Figure~\ref{fig:isotherm0} but using Model 1 instead. Since there
    is no known adsorption relation in this case, the Langmuir
    isotherm \eqref{eq:isotherm} in solid is displayed for reference
    only. A total of four unstable cases were detected (crosses).}
  \label{fig:isotherm1}
\end{figure}

\begin{figure}
  \ifpdf
  \includegraphics[width = 301pt]{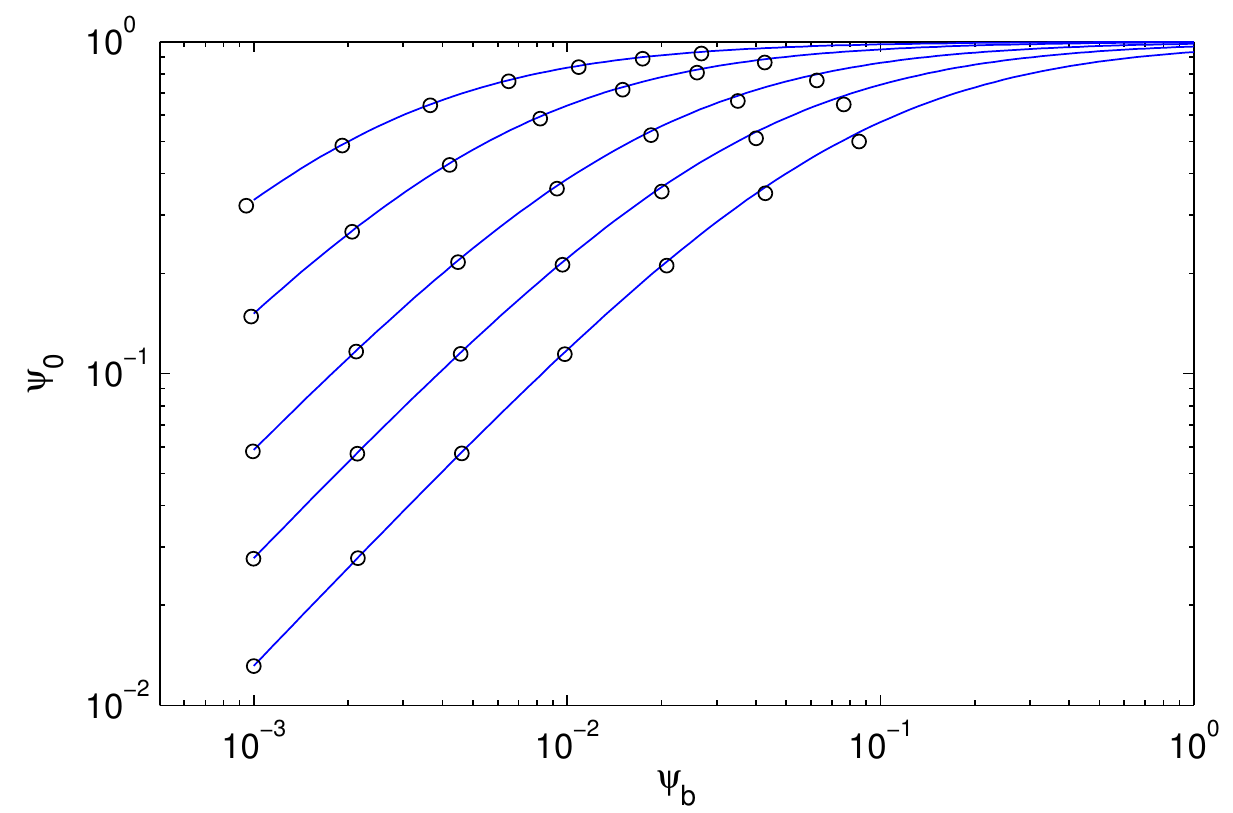}
  \includegraphics[width = 301pt]{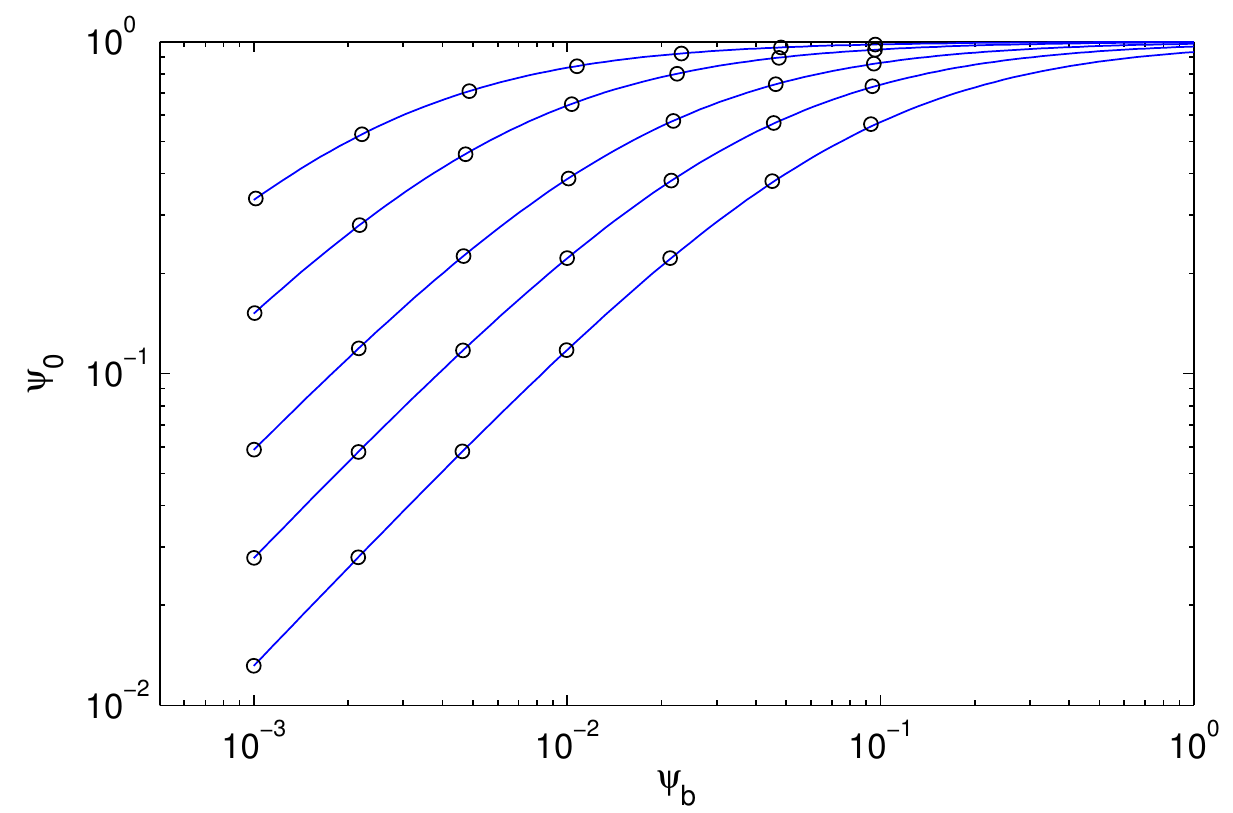}
  \else
  \includegraphics[width = 301pt]{fig/isotherm2.eps}
  \includegraphics[width = 301pt]{fig/isotherm3.eps}
  \fi
  \caption{Langmuir adsorption isotherm for Model 2 \textit{(top)} and
    3 \textit{(bottom)}.}
  \label{fig:isotherm2and3}
\end{figure}

\begin{figure}
  \ifpdf
  \subfigure[Model 1.]{\includegraphics[width = 103pt]{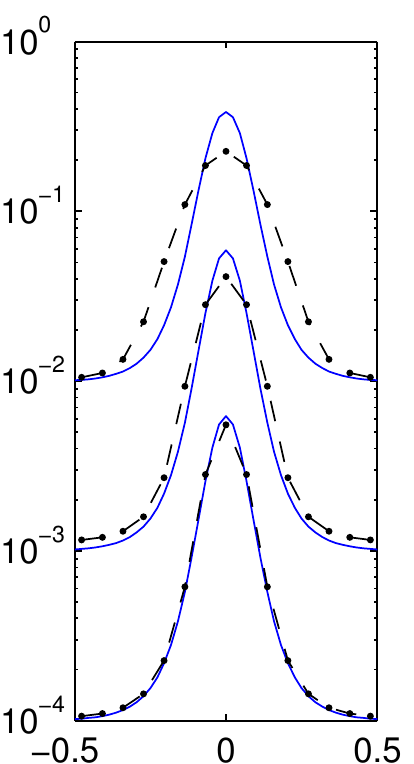}}
  \subfigure[Model 2.]{\includegraphics[width = 103pt]{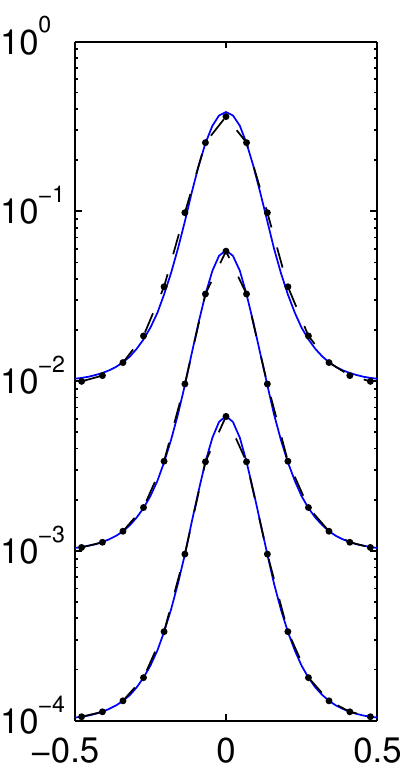}}
  \subfigure[Model 3.]{\includegraphics[width = 103pt]{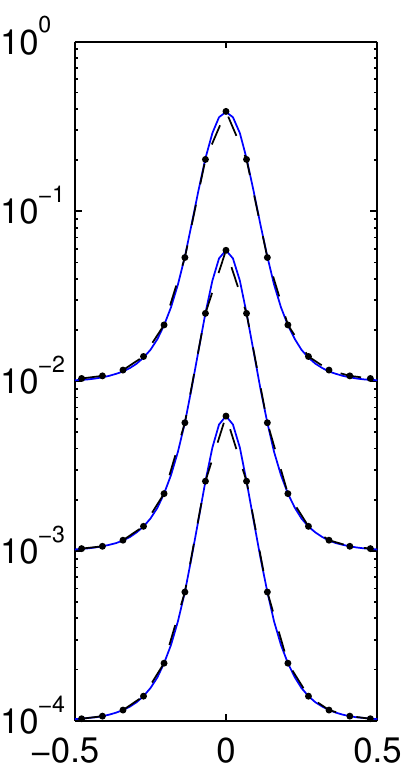}}
  \else
  \subfigure[Model 1.]{\includegraphics[width = 103pt]{fig/psisample1.eps}}
  \subfigure[Model 2.]{\includegraphics[width = 103pt]{fig/psisample2.eps}}
  \subfigure[Model 3.]{\includegraphics[width = 103pt]{fig/psisample3.eps}}
  \fi
  \caption{Initial surfactant concentration profiles according to
    \eqref{eq:psigen} \textit{(dashed)} and numerical equilibrium
    values \textit{(solid)} for different values $\psi_{b} \in
    \{10^{-4},10^{-3},10^{-2}\}$ and $\psi_{c} = 0.016$. For Model 1
    the value of $B-B_{b}$ in \eqref{eq:dB} formally belonging to
    Model 0 was used in order to obtain a closed expression for the
    initial profile.}
  \label{fig:samples}
\end{figure}

In Figure~\ref{fig:isotherm2and3} we similarly compare the theoretical
isotherm \eqref{eq:isotherm} with numerical values for Model 2 and
3. No unstable cases were detected and the results are all in very
good agreement. For Model 2, the initial data used in the experiments
were a bit off the actual equilibrium causing the measured values to
``creep'' slightly. However, all values stay in close proximity to the
predicted isotherm curve.

Finally, in Figure~\ref{fig:samples} the initial approximate
surfactant concentration profile is compared to the final equilibrium
profile for a single value of $\psi_{c}$ and for different surfactant
bulk concentrations $\psi_{b}$.

\begin{conclusion}
  We verified the theoretical analysis in
  Section~\ref{subsec:illposedness} and in particular the sharpness of
  the sufficient condition for ill-posedness \eqref{eq:psiblim}. The
  Langmuir isotherm \eqref{eq:isotherm} for Model 0 as published
  originally in \cite{diffsurf} and again, using altered parameters,
  in \cite{diffsurf2} could therefore not be obtained. Model 1,
  although stable in most cases, does not \textit{a priori} possess a
  natural adsorption isotherm. Model 2 and 3 are not only stable, but
  they also satisfy the Langmuir isotherm very accurately. Model 3
  produces the most sharp surfactant profile and Model 1 yields the
  most damped profile, at least towards the higher values of
  surfactant concentrations tested here. The analytical prediction
  \eqref{eq:psigen} is very accurate for Model 2 and 3.
\end{conclusion}

\subsection{Adsorption dynamics}

The adsorption dynamics for an interface in contact with a
semi-infinite bulk was considered in the early paper
\cite{ward_tordai} by Ward and Tordai. Here the time-dependent
decrease of surface tension under the presence of a solute was
explained on the basis of diffusion. Their approach has later been
refined; notably by Hansen \cite{hansen} who included also the process
of evaporation of solute from the interface, and by Ariel, Diamant,
and Andelman \cite{surfkinetics0,sharpsurf} who proposed a separate
treatment for ionic surfactants.

In order to cohere with this classical set-up we continue to use our
one-dimensional setting with an interface at $x = 0$. We treat the
phase-field variable $\phi$ as being time-independent in the arguments
below; the precise form is not critical as long as the scaling is such
that the interface is located within $|x| \lesssim \Cn$ and such that
$|\phi(x)| \sim 1$ for $|x| \gg \Cn$. As before we assume the constant
scalar bulk concentration $\psi_{b} = \psi(t,x \to \infty)$ and
further put $\psi_{0}(t) := \psi(t,x = 0)$, and also $\psi_{0,\eq} :=
\psi_{0}(t \to \infty)$.

The model proposed in \cite{sharpsurf} is ``sharp'' in the sense that
it is directly formulated as a \emph{discrete} model with a
characteristic length-scale $a$ (compartment size) on the order of a
single molecule. We may therefore refer to such a model as
``microscopic'' and it is clear that any continuous model will have
difficulties as the sharpness of the interface approaches this
length-scale. Indeed, for the diffuse interface model under current
consideration, the strict limit $\Cn \to 0$ does not make sense as the
limiting equilibrium surfactant concentration profile necessarily
becomes a constant bulk-value with a single point removed at the
interface itself (this follows from the isotherm relation
\eqref{eq:psigen}). However, the limiting behavior as $\Cn \to a \ll
1$ is still of interest since it makes a \emph{dynamic} comparison
with the microscopic model in \cite{sharpsurf} possible. Specifically,
we are interested in testing the adsorption dynamics at the interface
for which there are known discrete modeling approximations available.

In order to estimate the behavior as $\Cn \to a$, where we may think
of $a > 0$ as the length-scale of a single molecule, we write the
surfactant concentration in terms of \emph{inner} and \emph{outer}
variables,
\begin{align}
  \label{eq:ydef}
  \psi(t,x) &= \psiin(t,y) \qquad \mbox{ for } -1 \le y \le 1,
  \mbox{ with } y := x/\Cn, \\
  \psi(t,x) &= \psiout(t,x) \qquad \mbox{ for } |x| > \Cn.
\end{align}
Using the constant phase-field profile $\phi(x) = \tanh(x/\Cn)$ we can
write the PDE \eqref{eq:psiPDE} for $\psi$ as
\begin{align}
  \frac{\partial \psi}{\partial t} &= \frac{\Pii}{\Pe_{\psi}}
  \frac{\partial^{2} \psi}{\partial x^{2}}+
  \frac{1}{\Pe_{\psi}} \frac{\partial}{\partial x} \left[ \psi(1-\psi)
    \varphi'(x) \right], \\
  \intertext{in terms of which}
  \varphi(x) &= \left\{
    \begin{array}{lr}
      -\frac{\left( 1-\phi^{2} \right)^{2}}{4}+
      \frac{\phi^{2}}{4\Ex} & \mbox{(Model 0 and 3)} \\
      -\frac{\left( 1-\phi^{2}\right)^{2}}{4}-
      \frac{\Cn^{2}}{2}\psi''-\frac{\sigma}{2}\psi+
      \frac{\phi^{2}}{4\Ex} & \mbox{(Model 1)} \\
      \frac{\phi^{2}}{4}+
      \frac{\phi^{2}}{4\Ex} & \mbox{(Model 2)}
    \end{array} \right..
\end{align}
For $|x| \gg \Cn$ and in all four cases we readily get the bulk
diffusion equation,
\begin{align}
  \label{eq:outer}
  \frac{\partial \psiout}{\partial t} &= \frac{\Pii}{\Pe_{\psi}}
  \frac{\partial^{2} \psiout}{\partial x^{2}}.
\end{align}
Also, through the change of variables $y = x/\Cn$ expressed in
\eqref{eq:ydef}, we get the inner dynamics
\begin{align}
  \label{eq:inner}
  \frac{\partial \psiin}{\partial t} &= \frac{\Pii}{\Cn^{2}\Pe_{\psi}}
  \frac{\partial^{2} \psiin}{\partial y^{2}}+
  \frac{1}{\Cn^{2}\Pe_{\psi}} \frac{\partial}{\partial y}
  \left[ \psiin(1-\psiin) \varphi'(y) \right],
  \intertext{where for Model 1 the change of variables implies}
  \varphi(y) &= -\frac{\left( 1-\phi^{2}\right)^{2}}{4}-
  \frac{\psi''}{2}-\frac{\sigma}{2}\psi+
  \frac{\phi^{2}}{4\Ex}.
\end{align}
Assuming an asymptotic expansion of the form $\psiin(t,y) \sim \sum_{j
  \ge 0} \Cn^{j} \psi_{j,\Cn}(t,y)$ we have to leading order in $\Cn$
that the inner variable is in equilibrium with the outer variable via
the Dirichlet (matching) boundary condition $\psiin(t,\pm 1) =
\psi_{1}(t) := \psiout(t,\pm \Cn)$. As in
\cite{surfkinetics0,sharpsurf} we may interpret the new variable
$\psi_{1}$ as the concentration in the \emph{sub-surface layer} and we
have explicitly enforced an even symmetry $\psi(t,x) = \psi(t,-x)$ in
order to remain compatible with the modeling in those references.

We connect the two variables by requiring that the total volume is
conserved (assuming conservative outer boundary conditions),
\begin{align}
  \frac{dV}{dt} &= \int_{|x| \ge \Cn} 
  \frac{\partial \psiout}{\partial t} \, dx+
  \int_{-1}^{1} \frac{\partial \psiin}{\partial t} \, \Cn dy =
  -\frac{2\Pii}{\Pe_{\psi}}\psi_{1}'+
  \frac{d}{dt} \int_{-1}^{1} \psiin \, \Cn dy = 0.
\end{align}
For Model 2 and 3, with $\psiin$ in equilibrium the inner variable can
in principle be determined explicitly as a function of $\psi_{1}$
thanks to an isotherm-like relation. Model 1 does not possess such an
isotherm relation, but in this case a steady-state solution of the
inner dynamics \eqref{eq:inner} plays the same role. Although an
analytic expression for the resulting integral is lacking we argue
that under the present scaling, the near surface dynamics is
essentially of discrete character in the two variables $\psi_{0}$ and
$\psi_{1}$. As an example, using the trapezoidal rule we get the
approximative relation
\begin{align}
  \label{eq:closing}
  \frac{\partial \psi_{0}}{\partial t}+
  \frac{\partial \psi_{1}}{\partial t} &= 
  \frac{2\Pii}{\Pe_{\psi}\Cn}\psi_{1}'.
\end{align}

Interestingly, \eqref{eq:outer} and \eqref{eq:closing} above
correspond directly to Eqs.~(2.11) and (2.12) in the derivation in
\cite{surfkinetics0}. These equations can be solved analytically in
terms of an integral relation such that asymptotic estimates can be
obtained. More precisely, when starting with a uniform profile $\psi(t
= 0,x) = \psi_{b}$ we have the asymptotic ``footprint'' of diffusion
\cite{sharpsurf},
\begin{align}
  \label{eq:long}
  \frac{\psi_{0}(t)}{\psi_{0,\eq}} &\sim 
  1-\left( \frac{\tau_{0}}{t} \right)^{1/2}, \qquad t \to \infty. \\
  \intertext{There is also an initial transient phase with a
    $\sqrt{t}$-dependence \cite{hansen},}
  \label{eq:short}
  \psi_{0}(t) &\sim \mbox{const.}+\left( \frac{t}{\tau_{1}} \right)^{1/2}. \\
  \intertext{Finally, as pointed out in \cite{sharpsurf}, with a
    uniform initial profile there is also an ultra-short linear
    transient \emph{before} \eqref{eq:short} is valid,}
  \label{eq:ultrashort}
  \frac{\psi_{0}(t)}{\psi_{b}} &\sim 1+\mbox{const.} \times t.
\end{align}

\subsubsection{Numerical dynamics}

This quite interesting typical adsorption behavior is displayed in
Figure~\ref{fig:scales}, again using the numerical method outlined in
Appendix~\ref{app:num1}. The initial very short linear transient is
the phase where the region just outside the interface quickly becomes
depleted of surfactant. In the following $\sqrt{t}$-dependent phase,
surfactant diffuses from outside the sub-surface region, enters it,
and then immediately adsorbs to the interface, keeping the
concentration in the sub-surface region in equilibrium. The asymptotic
regime occurs towards the end of this process when the interface
gradually becomes saturated and the adsorption therefore slows down.

\begin{figure}
  \ifpdf
  \subfigure[]{\includegraphics[width = 167pt]{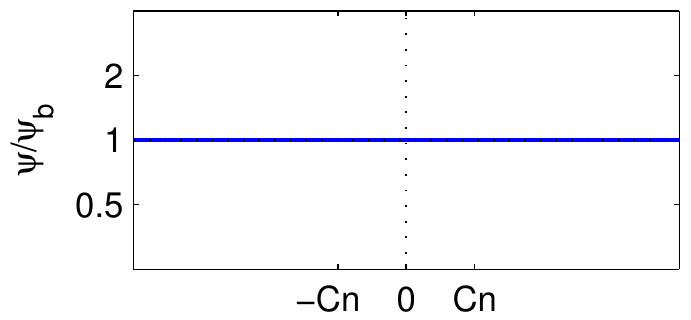}}
  \subfigure[]{\includegraphics[width = 167pt]{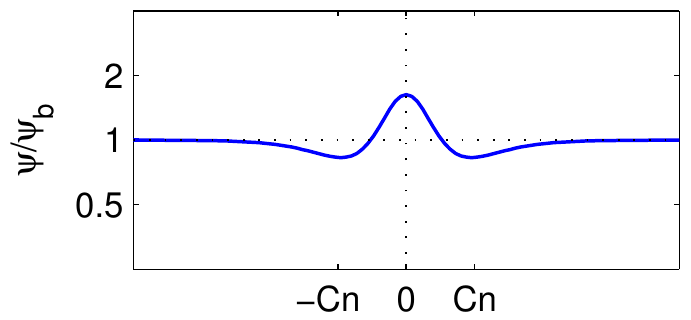}}

  \subfigure[]{\includegraphics[width = 167pt]{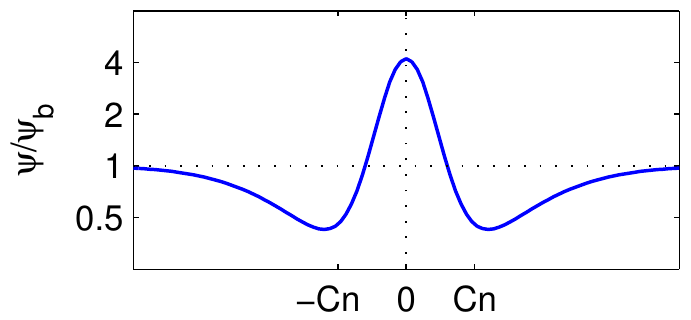}}
  \subfigure[]{\includegraphics[width = 167pt]{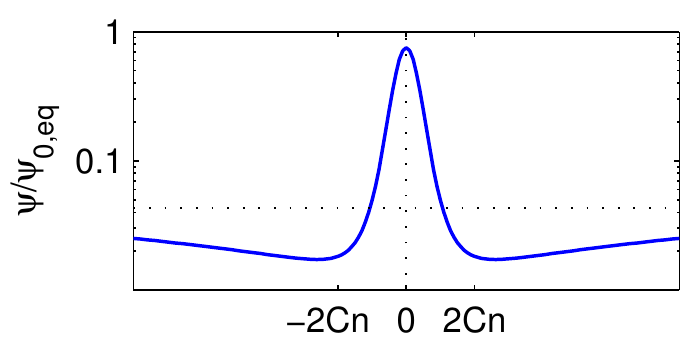}}

  \subfigure[]{\includegraphics[width = 169pt]{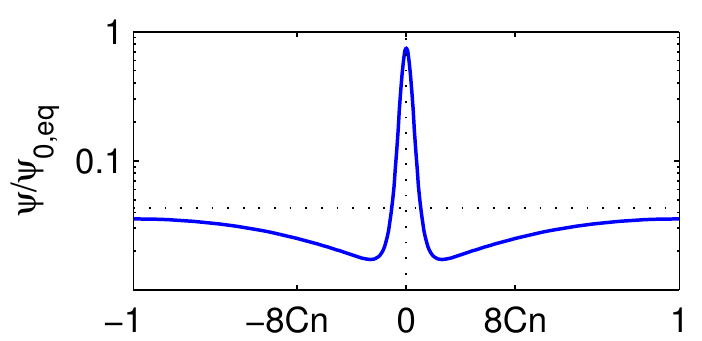}}
  \subfigure[]{\includegraphics[width = 169pt]{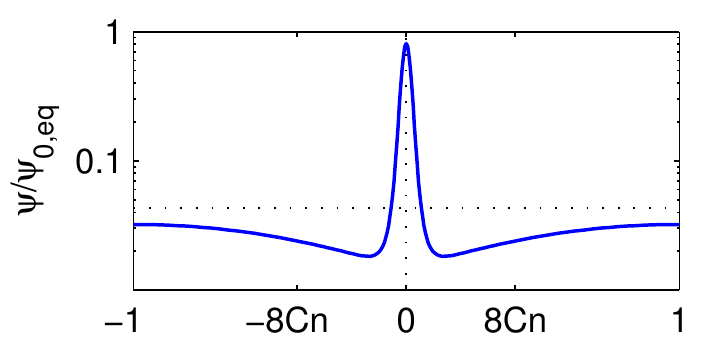}}
  \else
  \subfigure[]{\includegraphics[width = 167pt]{fig/ward_tordai_F31.eps}}
  \subfigure[]{\includegraphics[width = 167pt]{fig/ward_tordai_F32.eps}}

  \subfigure[]{\includegraphics[width = 167pt]{fig/ward_tordai_F33.eps}}
  \subfigure[]{\includegraphics[width = 167pt]{fig/ward_tordai_F34.eps}}

  \subfigure[]{\includegraphics[width = 169pt]{fig/ward_tordai_F35.eps}}
  \subfigure[]{\includegraphics[width = 169pt]{fig/ward_tordai_F36.eps}}
  \fi
  \caption{Adsorption phases. \textit{Top:} starting from a uniform
    concentration in (a), the interface very quickly drains the
    sub-surface region in (b)--(c). \textit{Middle:} surfactant now
    diffuses inwards from the bulk while the concentration at the
    interface steadily builds up in (c)--(d). \textit{Bottom:} in (e),
    the interface starts to reach saturation and in (f), finally,
    boundary effects due to the finite size of the numerical domain
    become prominent. In this example Model 3 was used with parameters
    $\Cn = 1/20$, $\psi_{c} = 0.016$, $\Ex = 1$, and $\psi_{b} =
    10^{-2}$. Note the logarithmic vertical scale.}
  \label{fig:scales}
\end{figure}

\begin{figure}
  \ifpdf
  \includegraphics[width = 297pt]{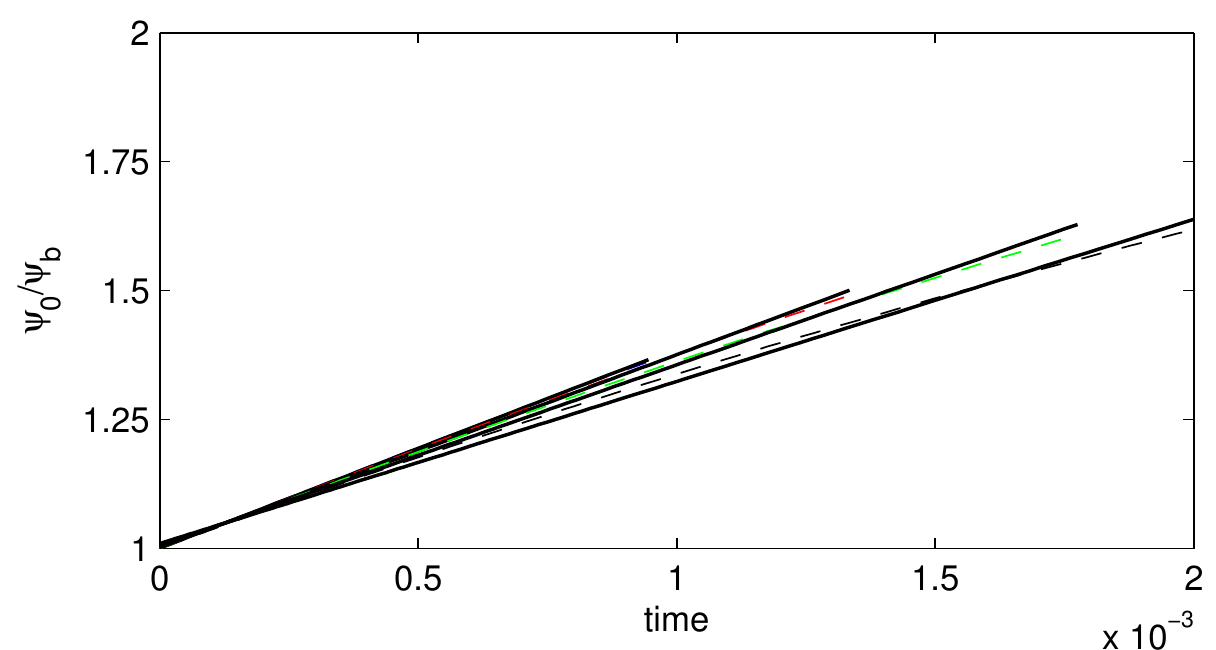}
  \includegraphics[width = 300pt]{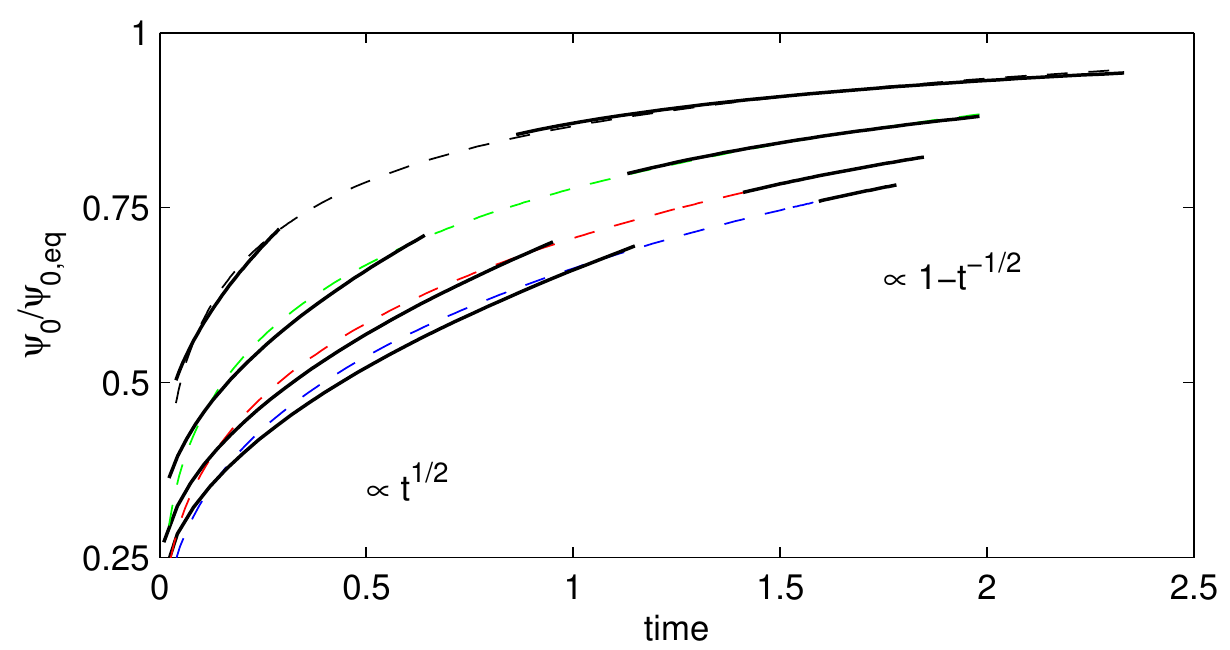}
  \else
  \includegraphics[width = 297pt]{fig/ward_tordaifitee_F2.eps}
  \includegraphics[width = 300pt]{fig/ward_tordaifit10_F2.eps}
  \fi
  \caption{\textit{Top:} initial very short linear transient, and
    \textit{bottom:} the characteristic diffusion-controlled
    $t^{1/2}$-transient together with the $t^{-1/2}$-asymptotics.
    \textit{Dashed:} numerical values, \textit{solid:} least-squares
    fit to the theoretical models (see text). These results are for
    Model 2 with parameters $\Cn = 1/20$, $\psi_{c} = 0.016$, $\Ex =
    1$. Counting from below in the bottom graph we have $\psi_{b} =
    [1,2,4,8] \times 10^{-2}$ (in the top graph the order is the
    reversed).}
  \label{fig:ward_tordai2}
\end{figure}

\begin{figure}
  \ifpdf
  \includegraphics[width = 296pt]{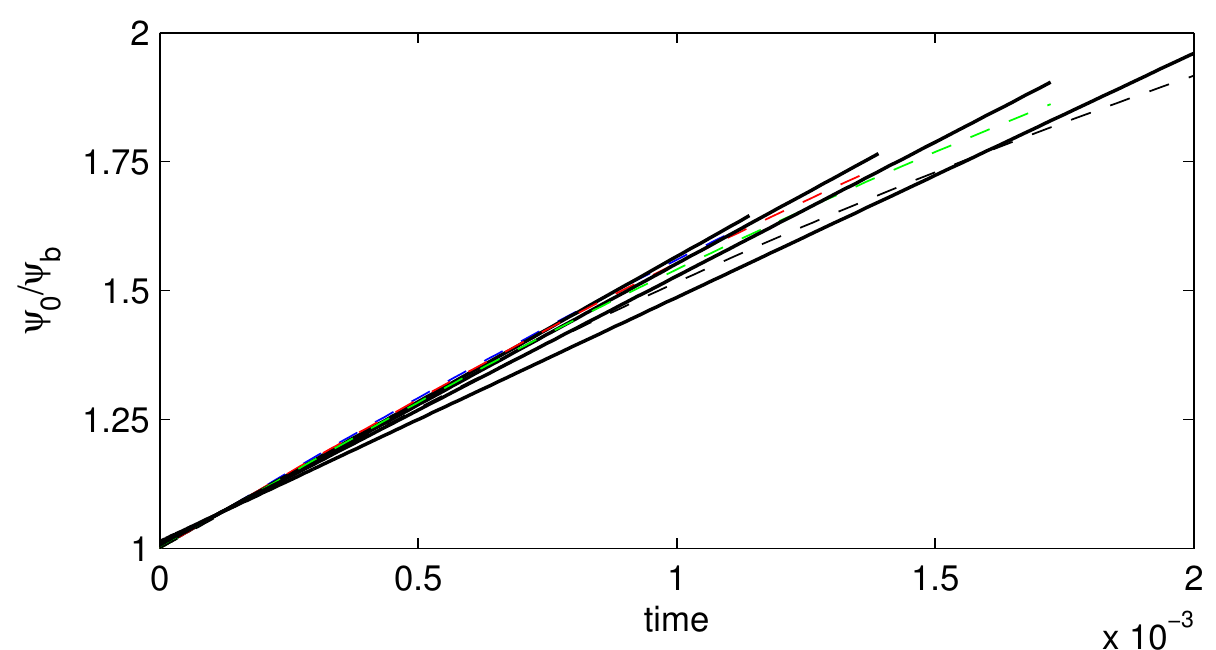}
  \includegraphics[width = 300pt]{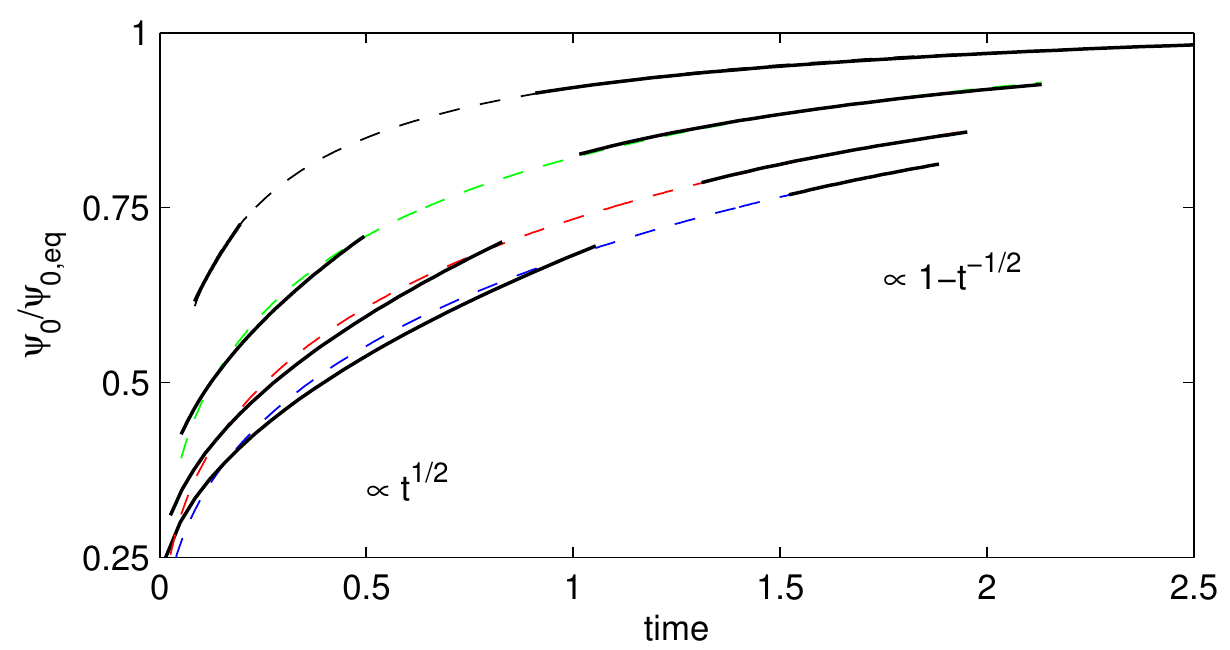}
  \else
  \includegraphics[width = 296pt]{fig/ward_tordaifitee_F3.eps}
  \includegraphics[width = 300pt]{fig/ward_tordaifit10_F3.eps}
  \fi
  \caption{As in Figure~\ref{fig:ward_tordai2} but for Model 3
    instead. Evidently, this model enjoys a slightly faster initial
    adsorption at the same Langmuir-controlled final saturation.}
  \label{fig:ward_tordai3}
\end{figure}

In Figure~\ref{fig:ward_tordai2} and \ref{fig:ward_tordai3} we
numerically fit (using standard polynomial least squares) the measured
time-dependent interfacial adsorption to the theoretical expressions
for a sharp interface expressed in
\eqref{eq:long}--\eqref{eq:ultrashort}. Although the two models
display some differences, both versions clearly fit quite well with
the theoretical predictions. Qualitatively similar results were
produced also for Model 1, but the temporal scaling is slightly
different here. Also, in one of the cases tested, we were unable to
obtain a solution due to instabilities and so we choose not to report
these results here.

\begin{conclusion}
  The Ward-Tordai governing equations \eqref{eq:outer} and
  \eqref{eq:closing} are approximately retrieved by performing a
  multiscale analysis in terms of inner and outer
  variables. Experimentally, for Model 2 and 3 we verified good
  agreement as to the theoretical behavior in the sharp interface
  limit.
\end{conclusion}


\section{Numerical method and experiments in 2D}
\label{sec:numexp2}

In this section we perform a final qualitative computational
experiment with the full hydrodynamic set of equations
\eqref{eq:phiPDE}--\eqref{eq:NS}. The purpose is to show that the
improved model can capture the nontrivial coupling between the fluid
flow and the surface forces. This is an essential property as it
involves the convective redistribution of surfactants, the resulting
gradients of surface tension that will appear, and a subsequent
modification of the overall dynamics. One very generic example which
we have chosen to study is to show that the presence of surfactant in
a high enough concentration may inhibit droplet coalescence.

We present the results from simulations in two spatial dimensions: the
extension to 3D would involve straightforward extensions of the
operators in the model equations. With the obvious exception of a
computationally heavier solution procedure, we do not foresee any
other particular differences, since the main physical effects of
surface stretching, convection of surfactant, diffusion, adsorption
and desorption are present also in the 2D test cases we have made.

Given the previous analytical and computational results we restrict
the experiment to Model 3 with chemical potentials given by
\eqref{eq:muphi3}--\eqref{eq:mupsi3}. Although we did observe the
ill-posedness for Model 0 here as well, those results are not detailed
here.

\subsection{Adaptive finite element method}

Successful finite element methods for the Cahn-Hilliard equations have
been devised by several authors. For a fully discrete convergent
method targeting the fully coupled hydrodynamic flow, see
\cite{CHNS_FEM}. Versions with degenerate mobilities are more
difficult to analyze such that \cite{degenerateCHFEM, numericalCHlog}
are two notable exceptions. Here, although only the no-flow case is
treated, on the other hand logarithmic terms are allowed in the free
energy.

We carried out our numerical simulations using femLego \cite{femLego},
a software to solve general partial differential equations with
adaptive finite element methods. The PDEs, the boundary conditions,
the initial data, and the method of solving each equation are all
specified in a Maple worksheet such that exact integration and
differentiation are possible.

The Cahn-Hilliard-type equations \eqref{eq:phiPDE}--\eqref{eq:psiPDE}
with chemical potentials according to
\eqref{eq:muphi3}--\eqref{eq:mupsi3} are treated as a coupled system
for the potentials $\mu_{\phi}$ and $\mu_{\psi}$ and the composition
variables $\phi$ and $\psi$. All equations are discretized in space
with piecewise linear functions and in time using the trapezoidal
rule, but with $\fatu$ frozen at the previous time-step, thus
resembling the strategy in \cite{CHNS_FEM}. The coupled system of
equations is solved using Newton iterations with UMFPACK
\cite{UMFPACK} as the inner linear solver.

To ensure mesh resolution in the vicinity of the interface, an
adaptively refined and derefined mesh is used with an \textit{ad hoc}
error criterion for each element $\Omega_{k}$,
\begin{align}
  \label{eq:errcrf}
  \int_{\Omega_{k}} \|\nabla \phi\|^{2} \le \TOL.
\end{align}
The mesh adaptivity is implemented by marking element $\Omega_{k}$ for
refinement if the element size is still larger than the minimum mesh
size allowed, $h > h_{\min}$, and it does not meet the criterion
\eqref{eq:errcrf}. In the case that an element does meet the
criterion, it is marked for derefinement unless it is an element of
the initial mesh. In the experiment reported here we adjusted $\TOL$
and $h_{\min}$ so that $\ge 10$ triangles were used across the
interfaces, defined here to be the regions where $\phi$ varies between
$-0.98$ to $0.98$. Some more details about this scheme can be found in
\cite{femLego}.

The Navier-Stokes equations \eqref{eq:cont}--\eqref{eq:NS} are solved
using the projection method devised in \cite{projectionFEM}. Again,
piecewise linear basis functions are used to discretize
space. Firstly, an approximate pressure is extrapolated from the
previous time-step and used when solving the momentum equation. Since
this by far is the most costly step, we employed an iterative linear
solver for this part \cite[\textit{GMRES},
  Chap.~3.2]{Greenbaum}. Secondly, a projection step in the form of a
Poisson equation derived from \eqref{eq:cont} is solved to correct the
velocity.

\subsection{Law of Laplace}

In order to test our 2D code we first performed numerical experiments
for the pure Cahn-Hilliard part (that is, $\psi = 0$ in
\eqref{eq:phiPDE}--\eqref{eq:NS}). Specifically, for a single droplet
in equilibrium we measured the pressure jump for different Cahn
numbers and different mesh spacings. The Cahn number is the ratio
between the width of the diffuse interface and the characteristic
length scale in the flow, here being the droplet diameter $d$. We kept
$\Ca = 1$, $\Pe_{\phi} = 3 \times 10^{-3}$, and $\Reynold = 1$ fixed
throughout these computations. According to Young-Laplace's equation
this gives an analytical pressure difference $\Delta P =
8\sqrt{2}/3$. The numerical domain was taken to be a cube of size $2d
\times 2d \times 2d$ and an equidistant mesh was
used. Table~\ref{tab:lap} shows the relative error between the
computed and the analytical pressure prediction for different Cahn
numbers and mesh spacings $h$, after equilibrium has been
approximately reached. Overall we observe a good agreement between the
numerical and the analytical predictions. Although the error in the
pressure depends on the numerical resolution of the interface we do
get a fairly accurate solution even with a comparably thick interface.

\begin{table}[!htb]
  \centering
  \begin{tabular}{lrrrrrrr}
    \hline
    $\Cn$ & 0.015 & 0.04 & 0.04 & 0.06 & 0.06 & 0.08 & 0.08 \\
    \hline
    $h$   & 0.003 & 0.013 & 0.02 & 0.02 & 0.03 & 0.027 & 0.04 \\
    $P_{error}$ & 0.06\% & 0.6\% & 2.0\% & 0.6\% & 1.8\% & 0.9\% & 1.7\% \\
    \hline
  \end{tabular}
  \vskip1mm
  \caption{Deviation between the numerical and analytical pressure for
    different Cahn numbers and mesh resolutions. Here $h$ is the
    (uniform) mesh spacing and $P_{error}$ is the relative error in
    the pressure jump.}
  \label{tab:lap}
\end{table}

\subsection{Colliding droplets}

\begin{figure}[htb]
  \centering
  \ifpdf
  \includegraphics[width=0.55\textwidth]{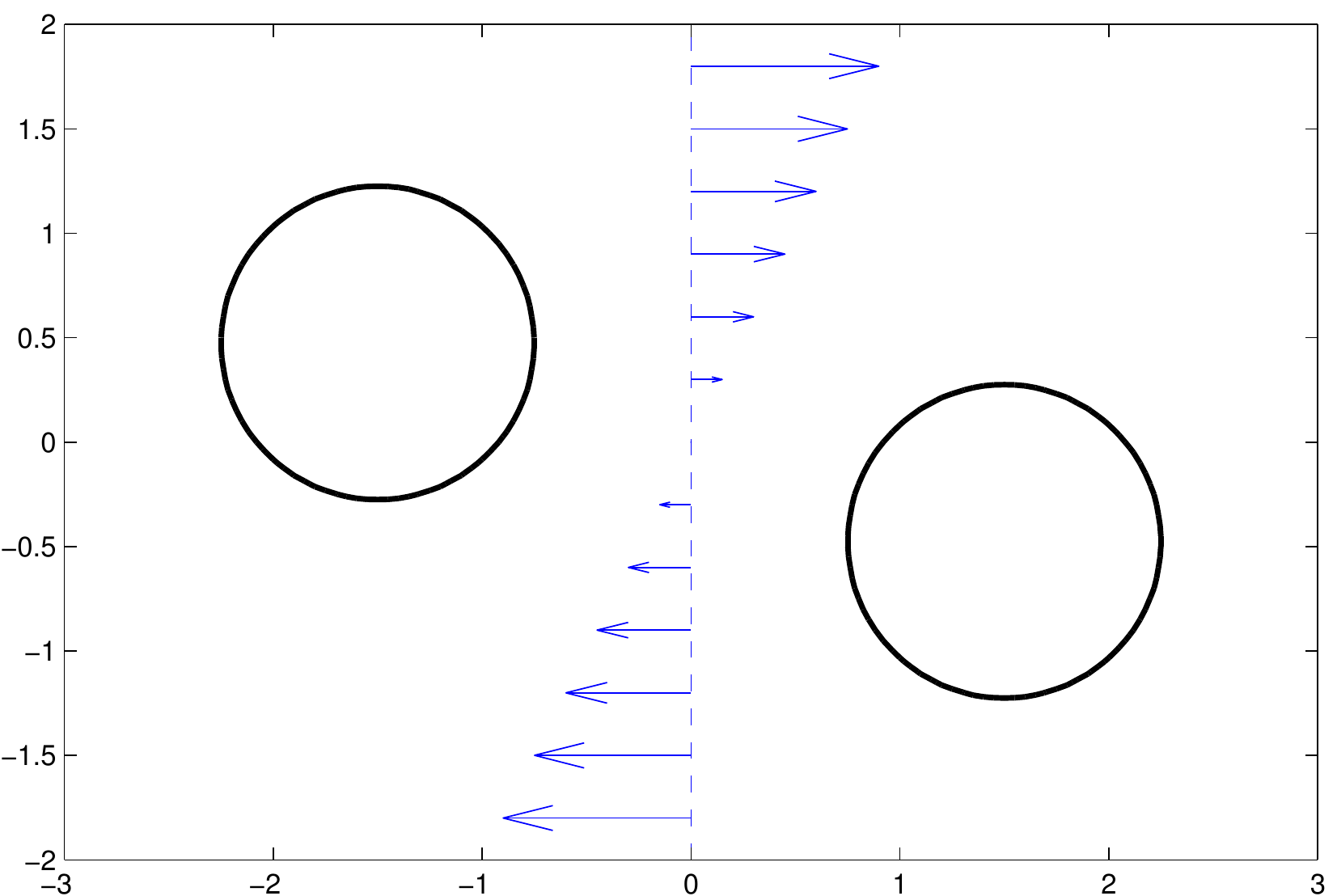}
  \includegraphics[width=0.55\textwidth]{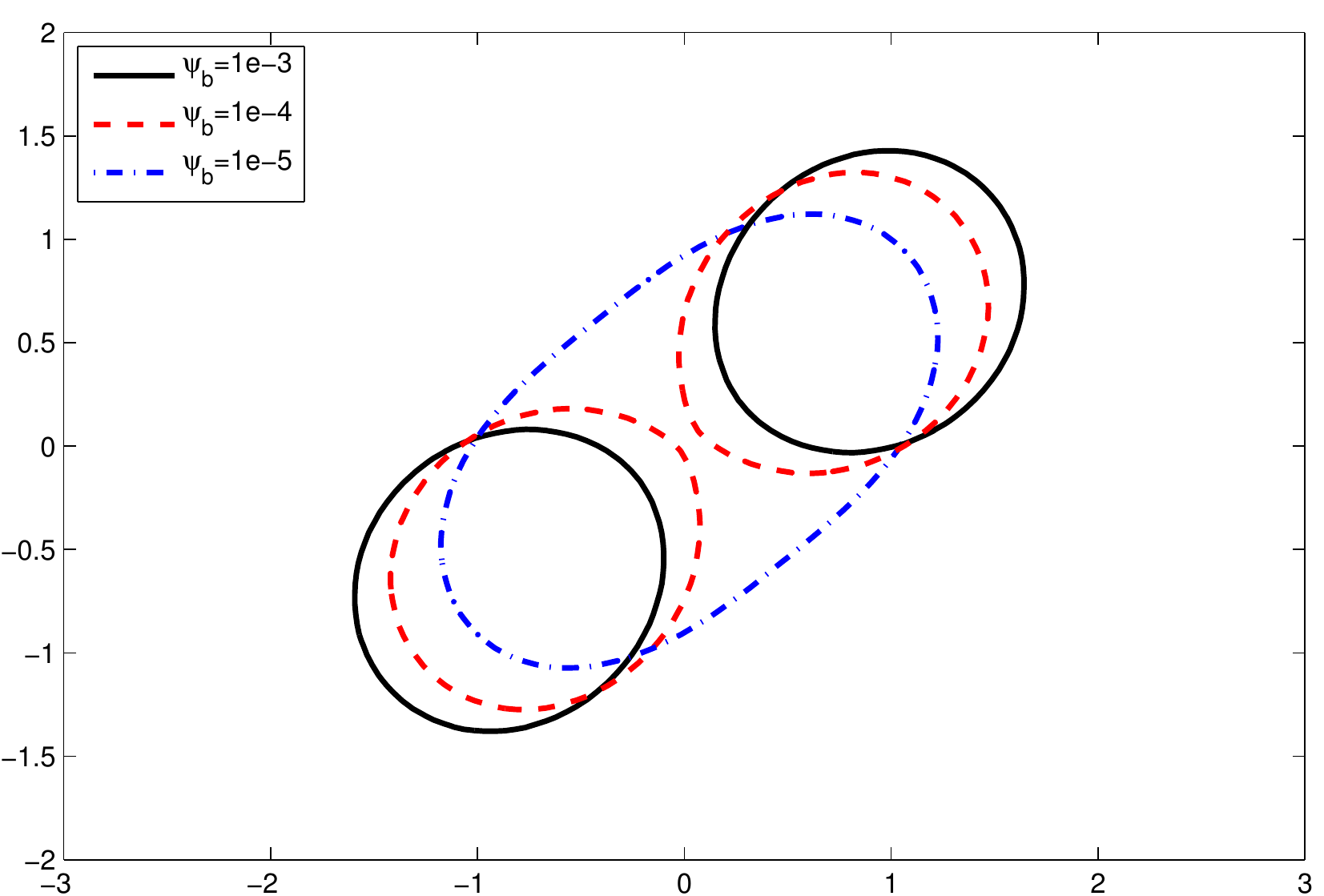}
  \includegraphics[width=0.55\textwidth]{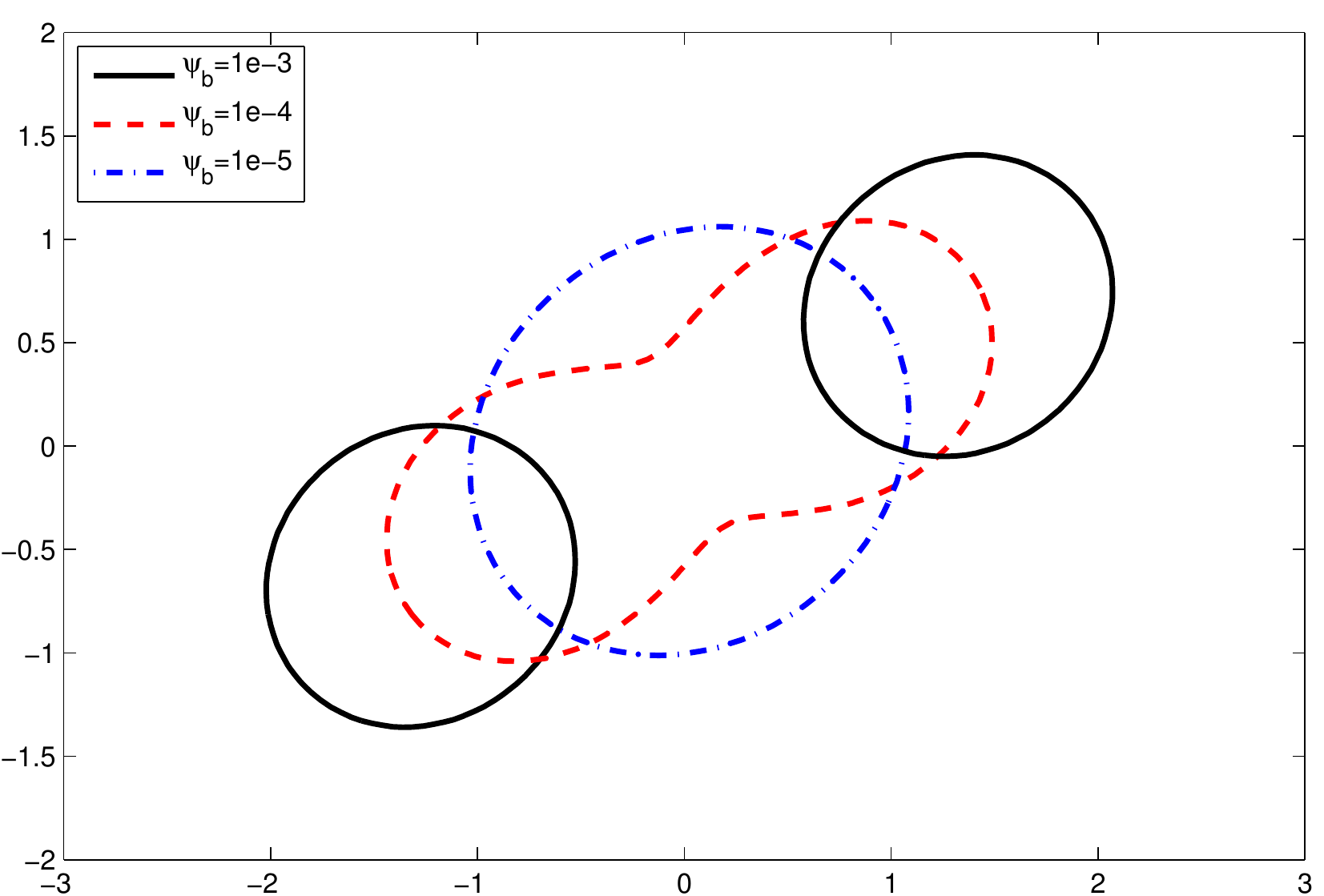}
  \else
  \includegraphics[width=0.55\textwidth]{fig/phi1_000.eps}
  \includegraphics[width=0.55\textwidth]{fig/phi1_165.eps}
  \includegraphics[width=0.55\textwidth]{fig/phi1_190.eps}
  \fi
  \caption{Evolution of two droplets in linear shear flow in a
    channel. Three cases with surfactant bulk concentration $\psi_{b}
    \in \{10^{-3},10^{-4},10^{-5}\}$ are displayed here. \textit{Top:}
    initial profiles; \textit{middle:} the solutions at time $t = 165$;
    \textit{bottom:} the solutions at $t = 190$. See text for further
    details.}
  \label{fig:phi1}
\end{figure}

\begin{figure}[htb]
  \centering
  \ifpdf
  \includegraphics[width=0.55\textwidth]{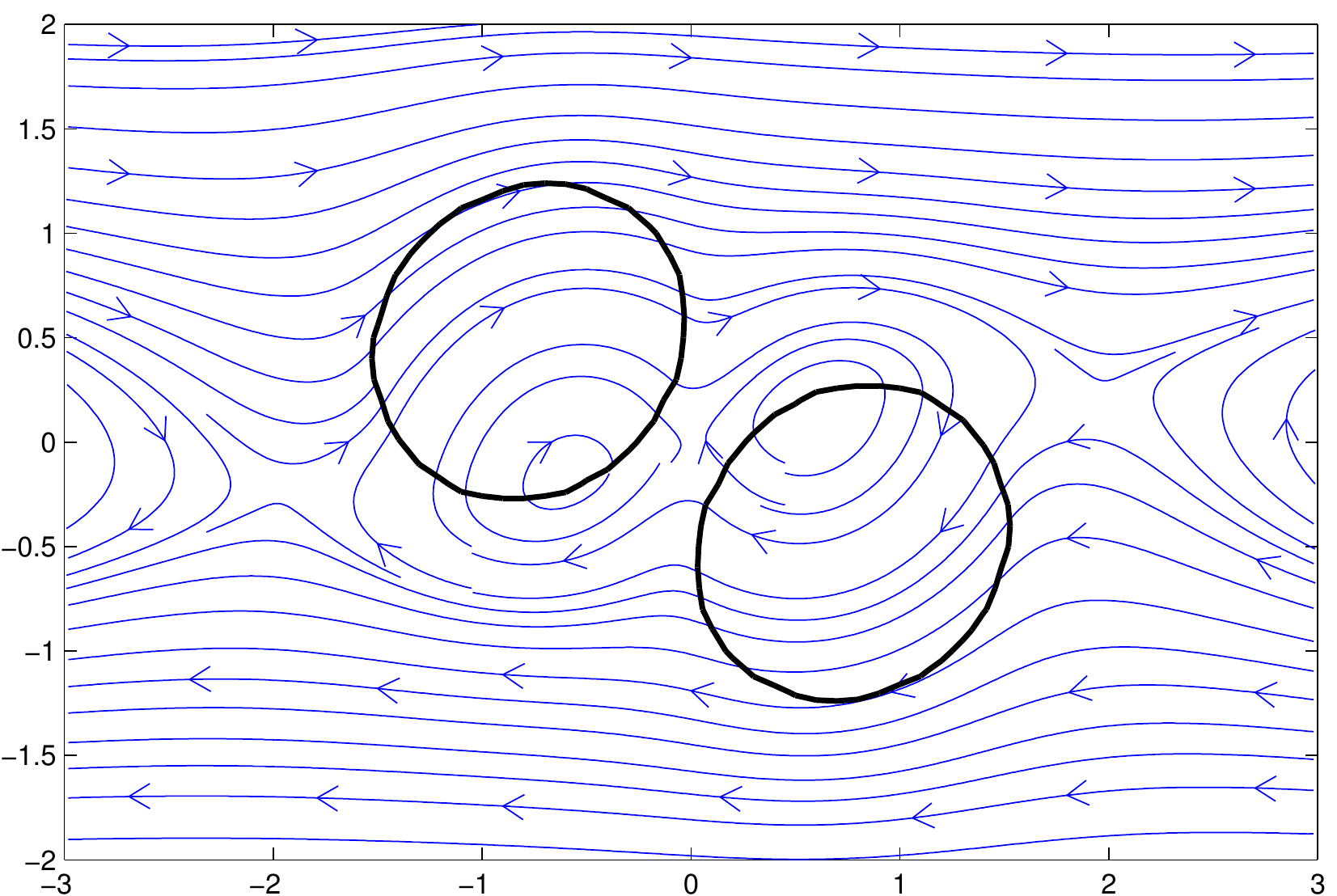}
  \else
  \includegraphics[width=0.55\textwidth]{fig/velocity.eps}
  \fi
  \caption{Complex flow (streamlines) for the two approaching droplets
    at $t = 75$ ($\psi_{b} = 10^{-3}$).}
  \label{fig:vel1}
\end{figure}

To demonstrate that our extended model of surfactant adsorption can be
coupled to hydrodynamics, we have used Model 3 to perform simulations
of two droplets colliding in a linear shear flow
channel. Figure~\ref{fig:phi1} shows the evolution of the two droplets
with droplet Reynolds number $\Reynold = 0.5$, capillary number $\Ca =
0.1$, P\'eclet numbers $\Pe_{\phi} = 10$, $\Pe_{\psi} = 100$, and
surfactant bulk concentration $\psi_{b} \in
\{10^{-3},10^{-4},10^{-5}\}$. The parameters of the surfactant model
for this simulation are given by $\Cn = 1/20$, $\Ex = 0.117$, and
$\Pii = 0.5857718$. All three cases are plotted on top of each other;
since increasing the amount of surfactant effectively hinders
coalescence of the droplets they are easily sorted out. For example,
with the smallest value of bulk concentration $\psi_{b} = 10^{-5}$,
the two droplets have coalesced already at $t = 165$. For $\psi_{b} =
10^{-4}$, they coalesce at $t = 190$ and they stay separate
indefinitely whenever $\psi_{b} \ge 10^{-3}$. Due to a comparably
strong surface tension effect and a weak fluid flow, we do not see a
noticeable effect on the droplet's behavior before they are in close
proximity to each other. However, when two droplets do approach, the
fluid flow becomes stronger and much more complex as shown in
Figure~\ref{fig:vel1}.

\begin{figure}[htb]
  \centering
  \ifpdf
  \includegraphics[width=0.55\textwidth]{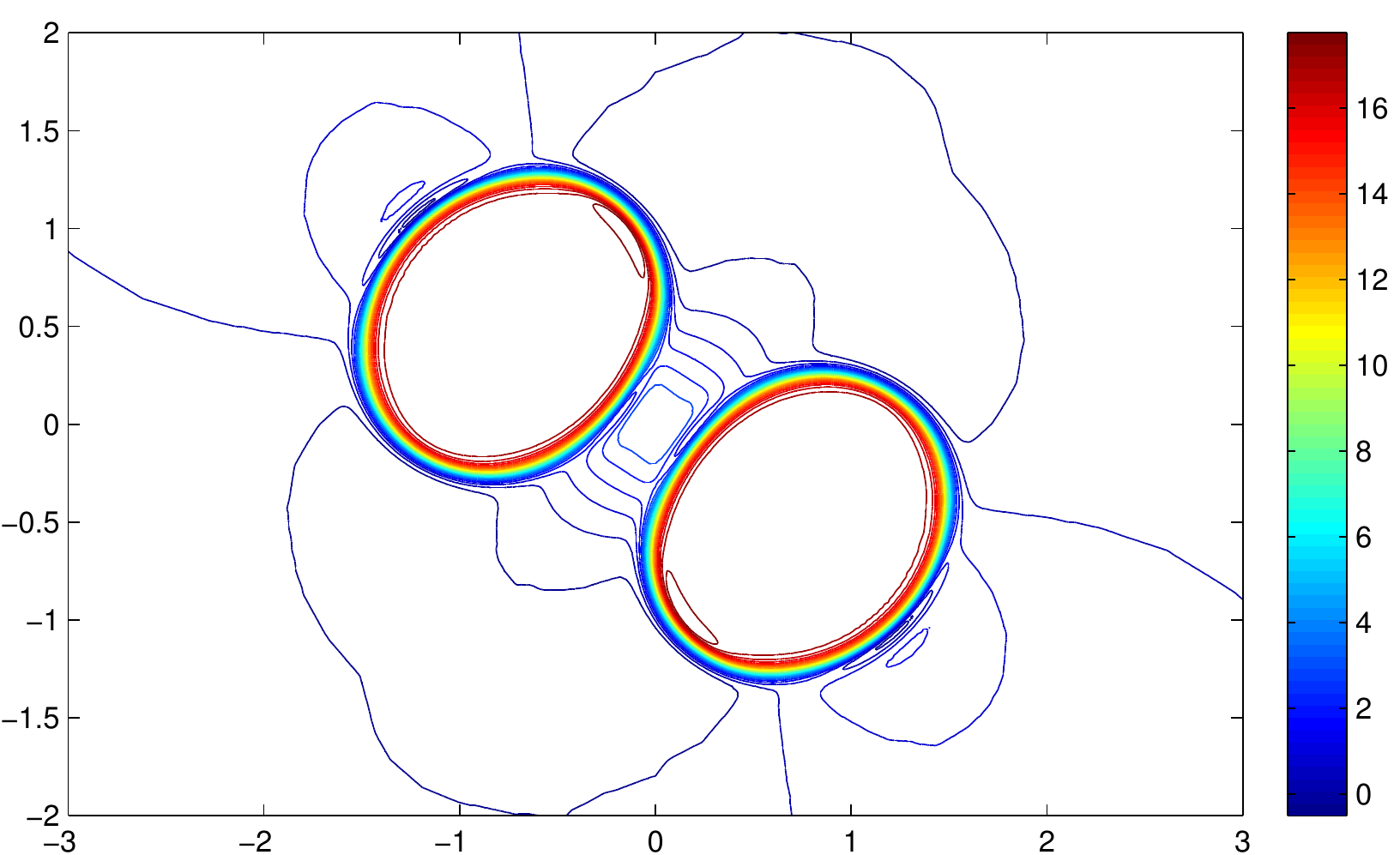}
  \else
  \includegraphics[width=0.55\textwidth]{fig/pressure.eps}
  \fi
  \caption{Pressure contour lines corresponding to the case in
    Figure~\ref{fig:vel1}.}
  \label{fig:pressure1}
\end{figure}

\begin{figure}[htb]
  \centering
  \includegraphics[width=\textwidth]{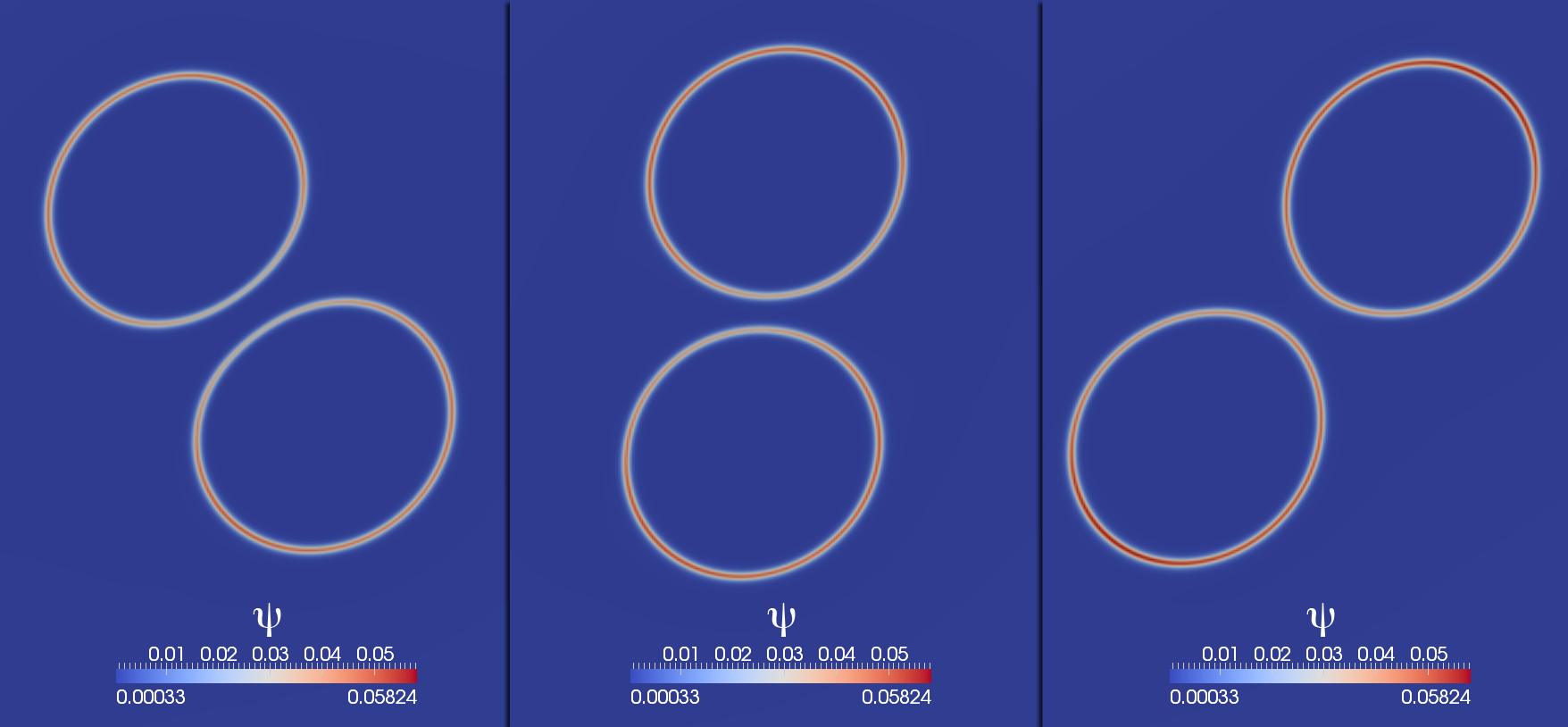}
  \caption{Surfactant concentration on the interface of the two
    droplets in linear shear flow. \textit{Left:} $t = 100$;
    \textit{middle:} $t = 125$; \textit{right:} $t = 150$ ($\psi_{b} =
    10^{-3}$).}
  \label{fig:phi2}
\end{figure}

Figure~\ref{fig:phi2} shows a contour plot of the surfactant
concentration of two colliding droplets in the simple shear flow with
$\psi_b = 10^{-3}$. Initially, the increased pressure in the gap
between the two droplets, see Figure~\ref{fig:pressure1}, pushes
surfactant away from the near-contact region, thus generating a
Marangoni stress that affects the droplet-droplet interaction. In
addition, the reduction of interfacial tension due to the presence of
surfactant has an effect on droplet deformation. In this way it
affects the droplet-droplet interaction also as a secondary effect.

\section{Summary and conclusions}
\label{sec:conclusions}

In \cite{diffsurf}, a diffuse interface model for the surfactant
adsorption onto the interface of two immiscible fluids was
presented. Recently, this model was extended to account also for
different solubility in the two phases, and the possibility to
consider systems better described by a Frumkin adsorption in addition
to the Langmuir adsorption \cite{diffsurf2}. However, the basic model
was the same.

Under the arguably quite weak assumption of a sufficiently smooth
equilibrium solution we have shown that this model is ill-posed in the
sense of frozen coefficients for a large set of physically relevant
parameters, in particular for cases where the interface becomes
saturated with surfactant. We have derived an explicit condition for
this ill-posedness, and through careful and very accurate
one-dimensional simulations illustrated that this condition is quite
sharp as indicated by instability and blow-up of numerical
solutions. These conclusions limit the usability of this model.

We have suggested and analyzed three alternatives to the basic model
(Model 0), as denoted Model 1--3. In Model 1, a natural idea in the
form of an energy contribution from the gradient of the surfactant
concentration is included. However, besides from being more
complicated, this model cannot reproduce the Langmuir adsorption
isotherm, and is also found to be numerically unstable in a few cases.
In Models 2--3, the surface energy accounting for the adsorption of
surfactant to the interface, is changed from a form containing the
gradient of the phase-field function, to a ``gradient-free form''.
This removes the problem of ill-posedness completely, and both models
are able to accurately reproduce the Langmuir adsorption
isotherm. Moreover, modeling different solubility in the two phases or
supporting the more general Frumkin isotherm would be straightforward
\cite{diffsurf2}. The exact form of the surface energy differs between
Models 2 and 3, rendering a somewhat different profile for the
surfactant concentration.

The conclusions concerning the mathematical models are also verified
in simulations with fluid flow. Interaction of droplets in two
dimensional shear flow in presence of surfactants is taken as a test
case. We find that the criteria for well-posedness determine the
usability of the model also in these more complex cases.

The phase-field modeling for surfactant laden flows is still in its
infancy. In this paper, we have highlighted the fact that mathematical
well-posedness of the resulting equations is not an obvious feature,
and that this is an issue that needs to be considered as new models
are developed.

\subsection{Reproducibility}
\label{subsec:reproduce}

Our 1D spectral-Galerkin code as described in Appendix~\ref{app:num1}
is available for download at the corresponding author's
web-page\footnote{\url{http://user.it.uu.se/~stefane/freeware}}. Along
with it, scripts that repeat the numerical experiments in
Section~\ref{sec:numexp1} and \ref{subsec:convergence} are
distributed.


\section*{Acknowledgment}

SE likes to acknowledge suggestions and comments by Per Lötstedt,
Gustaf Söderlind, Bertil Gustafsson, Sara Zahedi, and Jan Hesthaven.

This work was supported by the Swedish Research Council within the
FLOW and the UPMARC Linnaeus centers of Excellence. Computer time
provided by SNIC (Swedish National Infrastructure for Computing) is
gratefully acknowledged.

A-KT is a Royal Swedish Academy of Sciences Research Fellow supported
by a grant from the Knut and Alice Wallenberg Foundation and
thankfully acknowledges this support.

\subsection*{Author contributions}

SE wrote, developed the theory, and performed the experiments in
Section~\ref{sec:intro}, \ref{sec:analysis}--\ref{sec:numexp1}, and
Appendix~\ref{app:num1}. Section~\ref{sec:model} and
Appendix~\ref{app:nondim} were written by SE with inputs from
MD-Q. MD-Q wrote and performed the experiment in
Section~\ref{sec:numexp2}. A-KT and GA contributed to
Section~\ref{sec:intro} and wrote Section~\ref{sec:conclusions}.


\newcommand{\doi}[1]{\href{http://dx.doi.org/#1}{doi:#1}}
\newcommand{\available}[1]{Available at \url{#1}}
\newcommand{\availablet}[2]{Available at \href{#1}{#2}}

\bibliographystyle{abbrvnat}

\providecommand{\noopsort}[1]{} \providecommand{\doi}[1]{\texttt{doi:#1}}
  \providecommand{\available}[1]{Available at \texttt{#1}}
  \providecommand{\availablet}[2]{Available at \texttt{#2}}


\appendix

\section{Non-dimensionalization}
\label{app:nondim}

We discuss here in some detail the non-dimensionalization of the
diffuse phase-field model incorporating surfactants as originally
presented in \cite{diffsurf}. We follow the notation therein closely
in what follows.

Using the same degrees of freedom $\{\phi,\psi,\fatu\}$ as in
Section~\ref{sec:model} and with the same meaning, but now with $\phi
\in [-\phi_{0},\phi_{0}]$ we have (compare
\eqref{eq:phiPDE}--\eqref{eq:NS})
\begin{align}
  \label{eq:phiPDE0}
  \frac{\partial \phi}{\partial t}+\nabla \cdot (\phi \fatu) &=
  \nabla \cdot  M_{\phi} \nabla \mu_{\phi},  \\
  \label{eq:psiPDE0}
  \frac{\partial \psi}{\partial t}+\nabla \cdot (\psi \fatu) &= 
  \nabla \cdot M_{\psi} \nabla \mu_{\psi}, \\
  \label{eq:cont0}
  \nabla \cdot \fatu &= 0, \\
  \label{eq:NS0}
  \rho \left( \frac{\partial \fatu}{\partial t}+
  \fatu \cdot \nabla \fatu \right) &= -\nabla P+
  \nabla \cdot \left(\rho \nu [\nabla \fatu+(\nabla \fatu)^{T}] \right)
  -\phi\nabla\mu_{\phi}-\psi\nabla\mu_{\psi}.
\end{align}
Compared to \cite{diffsurf}, in \eqref{eq:cont0}--\eqref{eq:NS0} we
have changed the setting to incompressible flows for which $P$
enforces the incompressibility condition and the two final terms in
\eqref{eq:NS0} acts as surface tension forces \cite{CahnHilliard}. The
various free parameters of the model are defined as follows.

Firstly, the free energy still takes the form \eqref{eq:F} but now in
terms of
\begin{align}
  \label{eq:Fphi0}
  \Fphi &= -\frac{A}{2}\phi^{2}+\frac{B}{4}\phi^{4}+
  \frac{\kappa}{2}(\nabla \phi)^{2}, \\
  \Fpsi &= kT \left[ \psi\log \psi+(1-\psi)\log(1-\psi) \right], \\
  \Fone &= -\frac{\varepsilon}{2} \psi (\nabla \phi)^{2}, \\
  \Fex &= \frac{W}{2} \psi\phi^{2}.
\end{align}

In the scaling of \eqref{eq:Fphi0}, the equilibrium planar interface
at the origin is given by $\phi(x) = \phi_{0} \tanh (x/\zeta)$ with
$\phi_{0} = \sqrt{A/B}$, $\zeta = \sqrt{2\kappa/A}$, and with
$\{A,B,\kappa\}$ the parameters of the Cahn-Hilliard model
\cite{CahnHilliard}.

The chemical potentials are obtained through variational derivatives
of the free energy,
\begin{align}
  \mu_{\phi} &= \frac{\delta F}{\delta \phi} = 
  -A\phi+B\phi^{3}-\kappa\Delta\phi+
  \varepsilon\psi\Delta\phi+\varepsilon\nabla\psi\cdot\nabla\phi+
  W\psi\phi, \\
  \mu_{\psi} &= \frac{\delta F}{\delta \psi} = 
  kT \log \frac{\psi}{1-\psi}-
  \frac{\varepsilon}{2}(\nabla\phi)^{2}+\frac{W}{2}\phi^{2}.
\end{align}

To non-dimensionalize the model, let
$\{\hat{\phi},\hat{\psi},\hat{\fatu}\}(\hat{t},\hat{x})$ denote the
variables and coordinates as defined in the preceding paragraph. We
scale time and space in such a way that $(t,x) = (\hat{t}
u_{0}/L,\hat{x}/L)$ and change variables according to
$\hat{\phi}(\hat{t},\hat{x}) = \phi_{0}\phi(t,x)$,
$\hat{\psi}(\hat{t},\hat{x}) = \psi(t,x)$,
$\hat{\fatu}(\hat{t},\hat{x}) = u_{0} \fatu(t,x)$, and also set
$\hat{\rho} = \rho_{0} \rho$, $\hat{\nu} = \nu_{0} \nu$. It becomes
natural to define the new parameters
\begin{align}
  \label{eq:newparams}
  \Cn &= \zeta/L, \quad
  \Ex = \varepsilon/(W \zeta^{2}), \quad
  \Pii = kT/(A\phi_{0}^{2}).
\end{align}
Note that, as in \cite{diffsurf}, we tacitly assume that $\varepsilon
= \kappa$ in order to reduce unnecessary free parameters.

The new Ginzburg-Landau energy now becomes
\eqref{eq:F}--\eqref{eq:Fex}, yielding the chemical potentials in
\eqref{eq:muphi}--\eqref{eq:mupsi}. The new definition of the
mobilities is
\begin{align}
  M_{\phi} &= M_{\hat{\phi}} A, \quad 
  m_{\psi} = m_{\hat{\psi}} A \phi_{0}^{2}, \\
  \intertext{where we recall that a degenerate mobility $M_{\psi} =
    m_{\psi}\psi(1-\psi)$ is used for $\psi$ in \eqref{eq:psiPDE0}
    (conveniently, $M_{\psi} = \psi(1-\psi)$ in \eqref{eq:psiPDE}).
    It follows that the P\'eclet numbers in
    \eqref{eq:phiPDE}--\eqref{eq:psiPDE} are obtained as}
  \Pe_{\phi} &= L u_{0}/M_{\phi}, \quad \Pe_{\psi} = L u_{0}/m_{\psi},
  \intertext{and the Reynolds and Capillary numbers as usually by}
  \Reynold &= \frac{L u_{0}}{\nu_{0}}, \quad
  \Ca = \frac{\nu_{0}\rho_{0}u_{0}}{\gamma}.
\end{align}
Using the closing relation $\gamma \propto \kappa\phi_{0}^{2}/\zeta$
for the interfacial tension $\gamma$ finally yields \eqref{eq:NS}.


\section{Numerical method in 1D}
\label{app:num1}

In this section we describe our high-resolution one-dimensional
spectral scheme for the no-flow case ($\fatu = 0$ in
\eqref{eq:phiPDE}--\eqref{eq:psiPDE}), using the original energy
\eqref{eq:Fphi}--\eqref{eq:Fex}. The necessary modifications for Model
1--3 are trivial and are omitted for brevity.

A brief review of numerical methods for gradient flows in general and
surfactant models in particular might now be in order. Although many
papers treat numerical methods for the Cahn-Hilliard and the closely
related Allen-Cahn equations, it goes without saying that references
considering computational surfactant flows are much more scarce. In
this respect \cite{surfhybrid} stands out where a hybrid method for
surfactant two-phase flow is designed. In the setting of a sharp
interface treatment and using singular perturbation analysis, an
integral formulation is used in a semi-implicit time-discretization
strategy.

Semi-implicit, or split-step methods otherwise have a long history for
gradient flows, but is also an active area of research. Tracing their
origins back to general results in
\cite[\textit{(unpublished)}]{gradientstable} and
\cite{chgradientstable}, two more recent references include
\cite{biharmonicsplit} and \cite{chFEM}. In the former, general fourth
order gradient flows are considered, and in the latter the Allen-Cahn
and Cahn-Hilliard equations are specifically targeted.

A somewhat different design strategy is used in
\cite{allencahnhilliard,legendreCH} where the associated
Ginzburg-Landau functional is modified in order to provide for
numerical stability. The design therein relies on (spectral) Galerkin
formulations whereas in \cite{biharmonicsplit,CH3D} finite differences
are rather considered. A finite volume scheme for the Cahn-Hilliard
equation is obtained in \cite{timeCH}, where the efficiency benefit
with time-step adaptivity is also stressed.

On balance, and driven by a need for robustness and transparency
rather than for efficiency, we have chosen to design a polynomial
spectral Galerkin method with a very simple but carefully controlled
time-discretization strategy. We thus postpone more advanced spatial
and temporal simulation techniques for another occasion. As a word in
favor of this set-up, this discretization is similar to other
Galerkin-based methods such that generalizations to more realistic
situations are possible.

\subsection{Legendre spectral Galerkin method}

We consider semi-discrete weak solutions $\phi$, $\psi$ to
\eqref{eq:phiPDE}--\eqref{eq:psiPDE} ($\fatu = 0$) in the space
$\Hone(\Omega) \equiv \{w; \; \|w\|+\|\nabla w\| < \infty\}$ where
$\Omega = [-1,1]$ and the $\Ltwo(\Omega)$-inner product and induced
norm are understood.  Using natural (volume preserving) homogeneous
Neumann boundary conditions and the auxiliary variable $\Phi$ denoting
the chemical potential $\mu_{\phi}$ we get after integration by parts
the variational formulation
\begin{align}
  \label{eq:Phiweak}
  (\chi,\Phi) &= \left( \chi,-\phi+\phi^{3}+\frac{1}{2\Ex}\psi\phi
  \right)+\frac{\Cn^{2}}{2} (\nabla\chi,(1-\psi)\nabla \phi), \\
  \label{eq:phiweak}
  (\chi,\phi_{t}) &= -\frac{1}{\Pe_{\phi}}(\nabla\chi,\nabla\Phi)
\end{align}
for all $\chi \in \Hone(\Omega)$ and $t \in (0,T]$ assuming available
initial data at $t = 0$.

The variational formulation for $\psi$ is obtained in similar fashion
but using explicitly the observation in \eqref{eq:mupsiR}. Define
first
\begin{align}
  \label{eq:Psistrong}
  \Psi &= -\frac{\Cn^{2}}{4}(\nabla
    \phi)^{2}+\frac{1}{4\Ex}\phi^{2}.
\end{align}
We tacitly assume the boundary conditions $\Pii \partial \psi/\partial
n+\psi(1-\psi) \partial\Psi/\partial n = 0$ on $\partial \Omega$ and
get
\begin{align}
  \label{eq:psiweak}
  (\chi,\psi_{t}) &= -\frac{1}{\Pe_{\psi}}
  (\nabla\chi,\Pii\nabla\psi+\psi(1-\psi)\nabla\Psi),
\end{align}
again for all $\chi \in H^{1}(\Omega)$ and $t \in (0,T]$. In one
dimension these boundary conditions simplify to $\phi''' = \phi' =
\psi' = 0$ at the endpoints $x = \pm 1$.

As a discrete subspace of $\Hone$ we consider the space of polynomials
of degree $\le N$ and use as test- and trial functions the Legendre
polynomials which are orthogonal in the $\Ltwo([-1,1])$-inner product
\cite[Chap.~18]{DLMF10}. Expanding $\phi_{N}(t,x) = \sum_{n=0}^{N}
\phihat_{n}(t)\Leg_{n}(x)$ and similarly for $\psi_{N}$ and the
auxiliary variable $\Phi_{N}$ we get from the variational form
\eqref{eq:Phiweak}--\eqref{eq:psiweak} the semi-discrete set of
equations
\begin{align}
  \label{eq:Phiweakd}
  M \Phihat(t) &= a, \\
  \label{eq:phiweakd}
  M \phihat'(t) &= \alpha, \\
  \label{eq:psiweakd}
  M \psihat'(t) &= \beta.
\end{align}
The mass-matrix $M$ is diagonal with entries $M_{ii} = \|P_{i}\|^{2} =
2/(2i+1)$ for $i \in \{0,\ldots,N\}$ and
\begin{align}
  a_{i} &= \left( \Leg_{i},-\phi_{N}+\phi_{N}^{3}+
    \frac{1}{2\Ex}\psi_{N}\phi_{N}
  \right)+\frac{\Cn^{2}}{2} (\Leg_{i}',(1-\psi_{N})\phi_{N}'), \\
  \alpha_{i} &= -\frac{1}{\Pe_{\phi}}(\Leg_{i}',\Phi_{N}'), \\
  \beta_{i} &= -\frac{1}{\Pe_{\psi}}
  (\Leg_{i}',\Pii\psi_{N}'+\psi_{N}(1-\psi_{N})\Psi'),
\end{align}
with $\Psi$ defined by \eqref{eq:Psistrong} (hence imposed strongly).
The derivatives are readily obtained by explicit differentiation of
the Legendre expansions \cite[Appendix~B.1.4]{hesthaven_spectral}.

To evaluate the inner products we use the associated Gauss-Legendre
quadrature of at least the same order as the scheme, often
considerably higher in order to avoid or at least mitigate aliasing
errors \cite[Chap.~6]{hesthaven_spectral}. It is known that, for
nonlinear operators, aliasing errors and Gibbs oscillations near sharp
gradients may drive the scheme unstable even for smooth problems
\cite{legendrefilter}. This problem is perhaps mainly associated with
shocks but is of relevance also here since phase-field models
naturally produce sharp gradients in the vicinity of the interface.

We have therefore implemented \emph{filtering} along the lines
presented in \cite[Chap.~9.2]{hesthaven_spectral}. To this end we
replace the Legendre expansion with the \emph{filtered expansion}
\begin{align}
  \phi_{N}(t,x) = \sum_{n=0}^{N} \vartheta(n/N) \phihat_{n}(t)\Leg_{n}(x),
\end{align}
and similarly for $\psi_{N}$ and $\Phi_{N}$. A suitable filter is the
exponential one defined by \cite[p.~164]{hesthaven_spectral}
\begin{align}
  \vartheta(\eta) &= \left\{ \begin{matrix}
    1 & \eta \le \eta_{c}, \\
    \exp \left( -\alpha \left( 
        \frac{\eta-\eta_{c}}{1-\eta_{c}}\right)^{p} \right) 
    & \eta > \eta_{c}, \end{matrix} \right.
\end{align}
with $[\eta_{c},p,\alpha]$ filter parameters. Following the
recommendations in \cite{legendrefilter} we take the order to be $p =
12$ and set $\alpha = -\log \epsilon \approx 36$ with $\epsilon$ the
double precision machine accuracy. We also specify $\eta_{c} = 0.25$
so that the 25\% low modes are not filtered at all. The filtered
method is then conveniently implemented by replacing the mass-matrix
in \eqref{eq:Phiweakd}--\eqref{eq:psiweakd} with $\tilde{M}$ defined
by $\tilde{M}_{ii} := \vartheta(i/N)^{-1}M_{ii}$.

\subsection{Adaptive discretization of time}

Since we are considering problems that either develop instabilities or
for which there are nearby ill-posed problems, rather than aiming for
efficiency and high order we simply use the backward Euler method for
discretizing \eqref{eq:Phiweakd}--\eqref{eq:psiweakd} in time:
\begin{align}
  \label{eq:Phiweakdd}
  \tilde{M} \Phihat^{k+1} &= a^{k+1}, \\
  \label{eq:phiweakdd}
  \tilde{M} (\phihat^{k+1}-\phihat^{k}) &= \Delta t \alpha^{k+1}, \\
  \label{eq:psiweakdd}
  \tilde{M} (\psihat^{k+1}-\psihat^{k}) &= \Delta t \beta^{k+1},
\end{align}
with $\Delta t$ a suitable time-step and where the dependence on time
is expressed through superscripts.

Although there are arguably more efficient split-step methods
available that would probably apply here
\cite{biharmonicsplit,CH3D,chFEM}, the backward Euler method has a
specific advantage for gradient flows such as the one at hand. Namely,
it is one of very few \emph{unconditionally gradient stable} methods
(if not the only one), which preserves decay in the \emph{exact}
energy functional of the system (cf.~\cite[Chap.~5.6.1]{DSNA};
Theorems~5.6.1--5.6.3). Correct energy dissipation has been recognized
as one of the most important criteria when designing numerical methods
targeting gradient flows, which in specific cases may display chaotic
behavior and ``weak turbulence'' (see \cite[Sect.~4]{ODEstabrev} and
\cite{timeCH}). An intuitive argument is the fact that the
Ginzburg-Landau energy is \emph{the} modeling step and controlling the
rate of energy dissipation is therefore more important than formal
order.

For the Cahn-Hilliard equation it has been demonstrated by example
that time-step adaptivity is of vital importance since the dynamics is
rich in scales and since rapid transients due to initial data and
boundary conditions occur naturally \cite{timeCH,CH3D}. We address
this by using error estimates from half-steps and a careful time-step
selection mechanism.

Modern adaptivity based on \emph{digital filtering} was used in our
code \cite{dc,compstability}. To discuss it and for brevity, let
$T_{\Delta t}$ be the forward-in-time map by the backward Euler method
as obtained by solving \eqref{eq:Phiweakdd}--\eqref{eq:psiweakdd}. We
put $u_{1} = T_{\Delta t} u_{0}$ and $u_{2} = T_{\Delta t/2}^{2}
u_{0}$ for some initial state $u_{0}$. By Richardson extrapolation
this gives us a point-wise error estimate
\begin{align}
  \err &= |u_{1}-u_{2}|/3. \\
  \intertext{For given relative and absolute tolerances $\Rtol$ and
    $\Atol$ we compare the error to the effective tolerance vector
    $\TOL := \max(\Rtol|u_{2}|,\Atol)$ through the \emph{control
      variable}}
  \label{eq:cerr}
  \cerr &:= L(\|\err/\TOL\|_{\infty}^{-1/2}),
\end{align}
(element-wise division), where $1/2$ is the order of the error
estimate and where $L$ is a denoising \emph{limiter} (see below). We
consistently used the fairly stringent tolerances $\Rtol = 10^{-6}
\times \mbox{diag} (M^{-1/2})$ (vector) and $\Atol = 10^{-8}$ in our
experiments.

The control objective is now to keep $\cerr = 1$. An updated step-size
is obtained from a previous time-step $\Delta t'$ via the relation
$\Delta t = \rho \Delta t'$, where in terms of
\begin{align}
  \label{eq:rho}
  \rho &:= L(F(\cerr,\cerr',\rho'))
\end{align}
and where primes denote variables computed at the previous
time-step. As to the choice of $F$, the \emph{step-size controller},
we decided upon \textit{H211b} defined by
\begin{align}
  F(\alpha,\beta,\gamma) &= \alpha^{1/b}\beta^{1/b}\gamma^{1/b}
\end{align}
with the recommended choice of the single parameter being $b = 1/4$
\cite{dc}. Both \eqref{eq:cerr} and \eqref{eq:rho} employ the
\emph{arctan-limiter} \cite[Sect.~6]{compstability}
\begin{align}
  L(\alpha) &= 1+\kappa \arctan
  \left( \frac{\alpha-1}{\kappa} \right),
\end{align}
with, respectively, $\kappa = 2$ in \eqref{eq:cerr} and $\kappa = 1$
in \eqref{eq:rho}. The limiter's purpose in the latter case is to
smoothly restrict the rate by which the step-size may change and, in
the former case, to reduce the impact of noise in the error estimate.

By the very action of the step-size controller, the sequence of
$\rho$s is generally much smoother than the corresponding values for
the control variable in \eqref{eq:cerr}. For computational stability
and following a suggestion in \cite{compstability}, we therefore base
the decision of rejecting steps on the former variable rather than on
the latter. Specifically, given some tolerance factor $M$ (we took $M
= 2$), the controller is scheduled to reject the step whenever $\err >
M \cdot TOL$. Following the steps that lead to \eqref{eq:rho} this
suggests the \emph{rejection criterion}
\begin{align}
  \rho &< \rho_{\min} := L(F(L(M^{-1/2}),1,1)) \approx 0.9179.
\end{align}
In practice and with this strategy for time-step selection, step
rejections occur very rarely.

For solving the nonlinear system of equations
\eqref{eq:Phiweakdd}--\eqref{eq:psiweakdd}, simplified Newton
iterations with a numerical Jacobian were used. We followed the
prescription for estimating the speed of convergence in
\cite[Chap.~IV.8]{HNW2}, but also found it rewarding to augment the
code with control strategies devised in \cite{nlinODE} for reusing
previous Jacobian evaluations and factorizations.

The above described code was used successfully over a wide range of
test-cases. For reproducibility, the Matlab source code is made freely
available from the corresponding author's web-page. In particular,
scripts generating \emph{all} numerical results in
Section~\ref{sec:numexp1} are distributed (see
Section~\ref{subsec:reproduce} for further information).

\subsection{Convergence}
\label{subsec:convergence}

In order to assess the convergence properties of the code we performed
simulations for the pure Cahn-Hilliard part with $\Pe_{\phi} = 1$,
$\Cn = 0.1$, $u = 0.25$, and time $t \in [-0.5,0.5]$. In this case,
due to the non-zero flow $u$, there is an additional term
$(\chi,u\nabla \phi)$ in \eqref{eq:phiweak}. The reason for adding
transport is that for the corresponding free-space formulation there
is a traveling wave solution $\phi(t,x) = \tanh(x-u t)$, which in the
stated interval of time is an accurate solution also for the present
case with homogeneous Neumann boundary conditions. We parametrize the
tolerances for the time discretization according to $\Rtol =
\varepsilon \times \mbox{diag} (M^{-1/2})$ and $\Atol = \varepsilon
\times 10^{-2}$, with $\varepsilon$ a single parameter.

To test also the surfactant formulation
\eqref{eq:phiPDE}--\eqref{eq:psiPDE}, a highly accurate ($N = 400$,
$\varepsilon = 10^{-10}$) reference solution was computed for the
purpose of estimating errors. In this case a constant surfactant
profile $\psi = 10^{-3}$ was used as the initial data, and surfactant
parameters $\Pe_{\psi} = 1$, $\Ex = 1$, $\Pii = 0.1$ for Model 3 were
settled upon. The rest of the model was kept identical with the
previous case of pure Cahn-Hilliard flow.

\begin{figure}
  \ifpdf
  \includegraphics{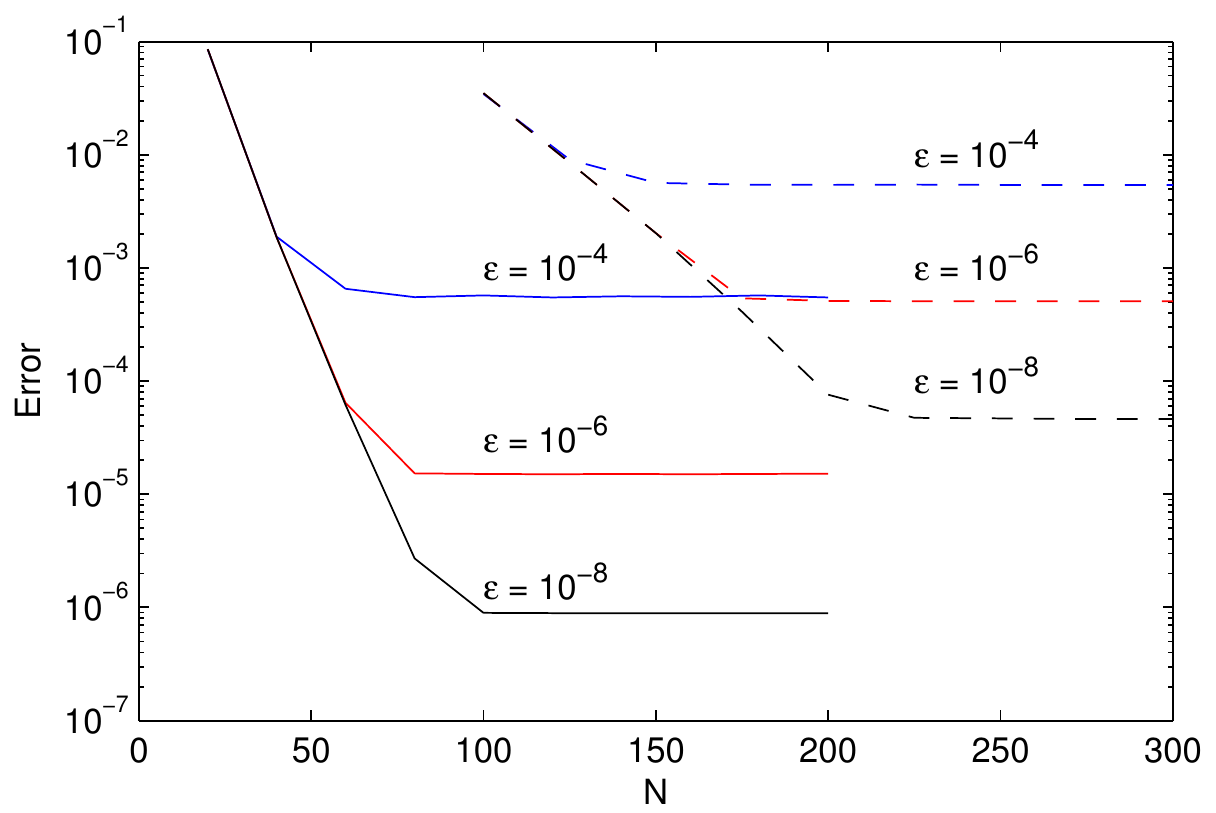}
  \else
  \includegraphics{fig/verification.eps}
  \fi
  \caption{Convergence study for the 1D code (relative error in the
    $L^{2}([-1,1])$-norm, maximum over the interval of
    integration). Three different time-stepping tolerances
    $\varepsilon$ and increasing order $N$ of the polynomial
    basis. \textit{Solid:} traveling wave solution (pure Cahn-Hilliard
    case), \textit{dashed:} coupled case, comparision with high order
    reference solution.}
  \label{fig:verification}
\end{figure}

The results of these experiments are shown in
Figure~\ref{fig:verification}. The coupled problem is evidently a more
difficult one, partially due to aliasing errors in this more strongly
nonlinear formulation, but also due to the relative sharpness of the
surfactant profile. In all cases spectral convergence down to a level
dictated by the accuracy of the time discretization are clearly
visible.

\end{document}